\documentclass[10pt]{elsarticle}
\usepackage[pdftex]{color}
\usepackage{xcolor}
\definecolor{darkyellow}{RGB}{200,200,0}

\usepackage [autostyle, english = american]{csquotes}

\usepackage[english]{babel} 
\usepackage{graphicx} 
\usepackage{amssymb}
\usepackage{amsmath}
\usepackage{caption}
\usepackage{subcaption}
\usepackage{bm}
\usepackage{algorithm}
\usepackage{algpseudocode}
\usepackage{todonotes}
\usepackage[pdfcreator = {pdflatex},
            hyperindex = {true},
            colorlinks = {true},
            linkcolor  = {blue},
            citecolor  = {blue},
            urlcolor   = {blue}]{hyperref}
\usepackage{cleveref}
\usepackage{multirow}
\allowdisplaybreaks
\newcommand{\bit}{\begin{itemize}}
\newcommand{\eit}{\end{itemize}}
\newcommand{\bma}{\begin{bmatrix}}
\newcommand{\ema}{\end{bmatrix}}
\newcommand{\m}{\mathbf}

\newcommand{\rmd}{\mathrm{d}}
\newcommand{\K}{\mathcal{K}}

\newcommand{\PS}[1]{\mathbb{PS}_{#1}} 
\newcommand{\OS}[1]{\mathbb{OS}_{#1}} 
\newcommand{\SL}{\mathcal{S}}
\newcommand{\SLt}{\mathcal{S}_\text{stk}}
\newcommand{\DL}{\mathcal{D}}
\newcommand{\SLP}[1]{\mathcal{S}_{#1}^{\mathbb{P}}}
\newcommand{\SLO}[1]{\mathcal{S}_{#1}^{\mathbb{O}}}
\newcommand{\DLP}[1]{\mathcal{D}_{#1}^{\mathbb{P}}}
\newcommand{\DLO}[1]{\mathcal{D}_{#1}^{\mathbb{O}}}
\newcommand{\SLPp}[1]{\mathcal{S}_{#1}^{\mathbb{P}\prime}}
\newcommand{\SLOp}[1]{\mathcal{S}_{#1}^{\mathbb{O}\prime}}

\newtheorem{theorem}{Theorem}

\newtheorem{definition}{Definition}

\begin{document}

\begin{frontmatter}
\title{Boundary integral equation analysis for spheroidal suspensions}
\author[1,*]{Leo Crowder}
\author[2,*]{Tianyue Li}
\author[1,**]{Eduardo Corona}
\author[2]{Shravan Veerapaneni}

\affiliation[1]{Department of Applied Mathematics, University of Colorado at Boulder}
\affiliation[2]{Department of Mathematics, University of Michigan}
\affiliation[*]{Equal contribution}
\affiliation[**]{Corresponding author; Eduardo.Corona@colorado.edu}
\date{}

\begin{abstract}
In this work, we provide a fast, spectrally accurate method for the evaluation of boundary integral operators (BIOs) on a suspension of prolate and oblate spheroids. 
We first derive formulas for the standard layer potential operators for the Laplace equation applied to an expansion of the integral densities in the appropriate spheroidal harmonic basis. These then lead to analytical expressions in solid harmonics that allow spectrally accurate evaluation of near-field particle interactions. Finally, a standard quadrature scheme is used to evaluate smooth, far-field interactions; these are then accelerated using the fast multipole method. 
Through a number of numerical test cases, we verify the accuracy and efficiency of our BIO evaluation framework for dense, polydisperse suspensions of spheroids. Through the use of standard formulas linking Stokes and Laplace potentials, we show our scheme can be readily applied to problems involving particulate suspension flows. For both Laplace and Stokes, our method allows us to evaluate BIOs for suspensions up to hundreds of particles on a single processor. 
\end{abstract}
\end{frontmatter}

\section{Introduction}
\label{sec:intro}

A wide class of soft materials can be effectively modeled as a discrete collection of mesoscopic particles, either by themselves (e.g. granular media) or suspended in and interacting through a physical medium (e.g. viscous fluid). The rich behavior of such materials stems from inter-particle interactions, both short-ranged and long-ranged, and can lead to complex emergent behavior at scales orders of magnitude larger than the particle size \cite{guazzelli2011physical, ness2022physics}. A key concern in the direct numerical simulation (DNS) of these systems is to resolve inter-particle interactions as well as other relevant physics in a way that does not compromise accuracy, stability, or favorable scaling of computational costs \cite{maxey2017simulation}. 

Inter-particle long-range forces are often well-captured by solving elliptic boundary value problems; this includes hydrodynamic forces for particles suspended in viscous flow (Stokes), and electromagnetic, hydrophobic and chemical reactions due to interparticle and particle-medium interaction (Laplace and screened Laplace). Numerical methods aimed at particulate system simulation can be classified with respect to their approach to long-range force fields; \emph{kernel smoothing methods} such as Rotne–Prager–Yamakawa \cite{wajnryb2013generalization}, \emph{approximation methods} such as Stokesian dynamics \cite{brady1988stokesian} and multipole methods \cite{cichocki1994friction}, and \emph{numerical PDE} methods such as immersed boundary \cite{uhlmann2005immersed}, fictitious domain \cite{patankar2000new}, and boundary integral methods (BIMs) \cite{pozrikidis}. 

In the absence of body forces or inhomogeneities in the surrounding medium, formulating these long-range force fields in terms of layer potentials involving boundary integral operators (BIOs)  for the associated PDE can be advantageous, as we need only solve a well-conditioned boundary integral equation at particle surfaces; this field and quantities of interest can then be computed to high precision using BIOs. For example, the single layer potential is defined by 
\begin{equation}
 \mathcal{S}[\sigma](\bm{x}) = \int_{\Gamma} G(\bm{x},\bm{y}) \sigma(\bm{y}) dS(\bm{y}), \label{BIO}  
\end{equation}
where $\Gamma$ is the particle boundary, $G$ is a free-space Green's function and $\sigma$ is a density function. Once discretized, 
fast summation methods such as the Fast Multipole Method (FMM) \cite{fmm3dbie} allow us to evaluate these operators with cost scaling linearly with the number of particles. This arguably makes BIM a natural fit for particulate media simulation, as demonstrated in their application to large-scale simulation frameworks for a wide range of particle shapes and behaviors in three dimensions \cite{rahimian2010petascale,lu2018parallel,yan2019scalable,malhotra2024efficient}.  


Spherical particle suspensions frequently serve as mimetic models for understanding the dynamics and rheology of a variety of complex fluids. Within active suspensions, for example, a mechanistic understanding of ciliated microorganisms such as {\em Opalina, Volvox,} and {\em Paramecium} was initially established through the {\em squirmer model}, which treats them as slip-driven spherical particles \cite{lighthill1952squirming, blake1971spherical, pedley2016spherical}. To account for shape anisotropy, the subsequent level of approximation often extends to spheroidal geometries; for example, the squirmer model was generalized to prolate spheroids \cite{keller1977porous}. These spheroidal models hold particular significance as they offer a valuable first approximation for various non-spherical microswimmers, including motile bacteria like {\em E. coli}. Naturally, spheroidal shapes have been widely employed in various applications to investigate the role of particle geometry on system behavior, such as in the design of optimal squirmers \cite{guo2021optimal, van2022effect}, the development of targeted drug carriers \cite{namdee2014vivo}, the analysis of diffusiophoretic particles \cite{poehnl2020axisymmetric} and the engineering of medical microrobots \cite{yoo2018actively, Kim_Jeon_Yang_Jin_Kim_Oh_Rah_Choi_2023}. For a comprehensive review of the differences between spherical and spheroidal particle systems in bioscience applications, see \cite{gadzinowski2021spherical}. The numerical algorithms presented here are broadly applicable to this wide range of applications involving spheroidal particles suspended in Stokes flow, spanning diverse scenarios such as static particles (as in porous media), force- and torque-driven particles (as in mobility problem), slip-driven particles (as in active suspensions), and imposed velocity-driven particles (as in resistance problem) \cite{pozrikidis}.

\paragraph{Related work.} 
Here, we only highlight the prominent high-order approaches one can employ for evaluating \eqref{BIO} on spheroidal particles. The most general approach is adaptive triangulation of the boundary surfaces \cite{greengard2021fast}, followed by high-order singular quadrature \cite{bremer2012nystrom}, Quadrature by Expansion (QBX) \cite{klockner2013quadrature} or adaptive refinement, as implemented in the open-source software \texttt{fmm3dbie} \cite{fmm3dbie}. The recently developed recursive reduction quadrature rules \cite{zhu2022high, jiang2024recursive} can further improve the efficiency of the near-singular integration schemes. However, for simple shapes suspended in Stokes flow, global surface representations such as spherical harmonics are often more efficient owing to their spectral accuracy. They have been employed in simulating deformable particle systems such as bubble, capsule or vesicle flows \cite{zhao2010spectral, veerapaneni2011fast, rahimian2015boundary, sorgentone2018highly}. While  quadrature based on spherical grid rotations achieves spectral accuracy in evaluating singular integrals \cite{gimbutas2013fast}, high-order close evaluation schemes for global representations have been explored in \cite{sorgentone2018highly, nitsche2025corrected}. The former self-interaction evaluation scheme has a computational cost that scales as $\mathcal{O}(m^2 \log m)$ for $m$ discretization points. Owing to the spectral accuracy, however, $m$ is much smaller compared to local patch, or triangulation, based discretization schemes for a target precision.  Since in our target application, the particle shapes are fixed, and moreover axisymmetric, BIOs can be reduced to convolutions with modal Green's functions on the generating curve. We can apply the numerical scheme introduced in \cite{young2012high} (for Laplace BIOs) to reduce the computational cost per particle to $\mathcal{O}(m^{1.5})$ via the use of FFTs in the azimuthal direction. 

The approach presented in this paper further reduces the complexity to $\mathcal{O}(m)$ for oblate and prolate spheroidal shapes by diagonalizing the self-interaction BIO (i.e., the BIO in \eqref{BIO} with $\bm{x} \in \Gamma$). 
It builds upon a similar approach developed in~\cite{vico2014boundary,corona2018boundary,kohl2023fast} for Laplace, Stokes, Helmholtz and Yukawa BIOs on spherical particles. One of the main advantages of this approach is that both for targets on and close to the surface $\Gamma$, the BIO \eqref{BIO} can be evaluated using {\em analytical formulas} with spectral accuracy, avoiding the need for nearly singular integral evaluation. 
Notably, this approach was leveraged in \cite{yan2019scalable} to simulate an active matter system involving around eighty thousand closely interacting spherical particles. The present work generalizes this approach to spheroidal particles.

We note that the {\em Method of Fundamental Solutions} (MFS) is another efficient approach for simulating Stokes flow around a collection of simple geometries. Recently, Broms et al.  \cite{broms2024method, broms2024accurate} significantly advanced this technique by accelerating MFS-based solvers for the Stokes mobility problem, enabling large-scale simulations. The essential idea is to avoid singular integrals in BIOs by placing proxy sources within the constituent rigid particles and enforcing boundary conditions at a set of collocation nodes. While this results in an ill-conditioned least squares system, its efficient, large-scale solution constitutes the core contribution of \cite{broms2024method}. Comparing MFS with the BIO analysis approach developed here reveals distinct advantages for each. While MFS is often simpler to implement, can handle arbitrary but simple shapes, and faster for low to moderate precision (e.g., up to six digits), the latter is beneficial for high-precision calculations owing to spectral accuracy, constrained geometry flows (as in Fig. \ref{fig:spheroids}) and its ability to use classical second-kind integral equation formulations \cite{hsiao2008boundary, corona2017integral} for various physical problems. In addition, we note that BIO analysis has several ancillary applications such as in constructing preconditioners (e.g., \cite{veerapaneni2009boundary, veerapaneni2011fast}), designing fast solvers (e.g., FMPS \cite{gimbutas2013fast}) and facilitating perturbative analysis such as the small deformation theory (e.g., \cite{vlahovska2009small, vlahovska2019electrohydrodynamics}).

\begin{figure}[ht!]
    \centering
    \includegraphics[width=\textwidth]{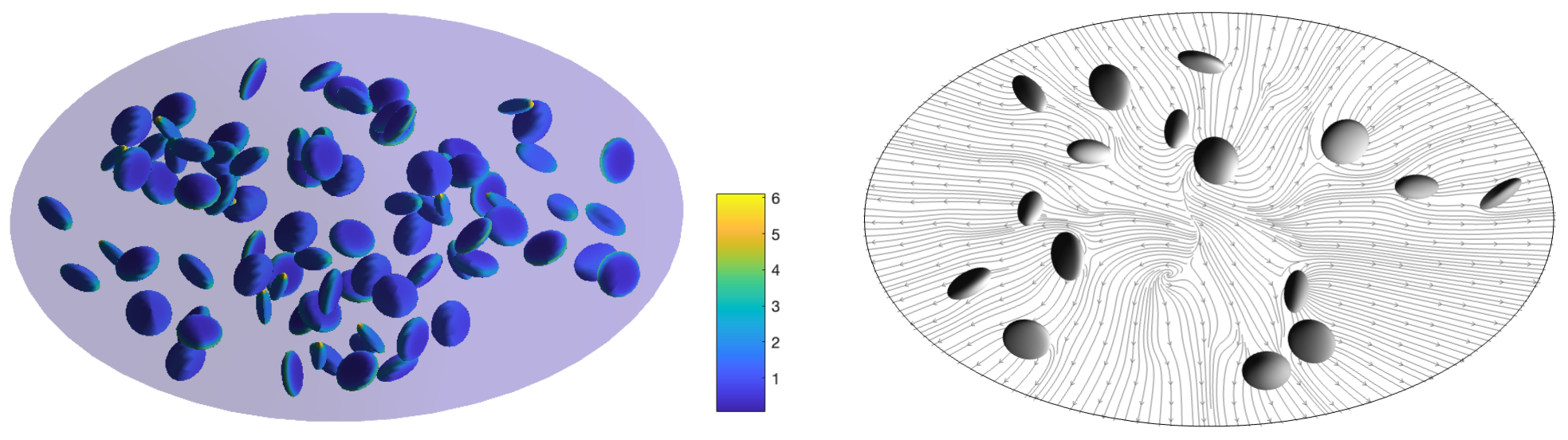}
    \caption{{\em \small Solution of the Stokes equations within a spheroidal confinement, driven by interfacial forces acting on a collection of both prolate and oblate spheroidal particles. (Left) Color indicates the magnitude of interfacial forces, here set to $\bm{f}^{\text{int}} = H\bm{n}$, where $H$ is the mean curvature and $\bm{n}$ is the unit normal on the surface of each spheroid. (Right) Plot of streamlines on a cross-section of the domain along with the particles intersecting this plane. The methods developed in this work allow us to compute close particle-particle hydrodynamic interactions as well as the flow field  arbitrarily close to spheroidal particle and geometry surfaces. }} \vspace{-0.1in}
    \label{fig:stokes flow}
\end{figure} 

\paragraph{Contributions.}

The methods contributed here extend the fast BIO evaluation approach in \cite{corona2018boundary,yan2019scalable} to layer potentials defined on a diverse collection of spheroidal prolates and oblates. 
Our key contributions include: 

\begin{itemize}
\item For an integral density $\sigma$ defined on a spheroidal surface, we find bases related to surface spheroidal harmonics such that the map from the expansion of $\sigma$ in this basis to that of Laplace layer potentials on and off the surface in spheroidal harmonics is \emph{diagonal}. 

\item Exterior, interior and on-surface evaluation formulas are derived for Laplace single and double layer potentials, as well as their first-order directional derivatives. Thereby, the task of evaluating weakly- and nearly-singular integrals is reduced to simple formula evaluation.   

\item We develop fast algorithms for both BIO matrix-vector product and surface-to-target potential evaluation for Laplace BIOs by splitting the task into near and far-field, handled by the spheroidal harmonic method and a smooth quadrature followed by the FMM, respectively. 

\item Using formulas that express Stokes layer potentials as a linear combination of Laplace BIOs \cite{tornberg2008fmm3dstokes}, we demonstrate the straightforward application of our methods to solving Stokes boundary value problems.

\item As part of our numerical studies, we contribute a heuristic scaling scheme aimed at minimizing conditioning of completed double-layer operators on moderate to high aspect ratio prolates and oblates. We find the slender body limit scaling in \cite{malhotra2024efficient} to be a close match for elongated prolates. 
\end{itemize}

\paragraph{Limitations.} 
By construction, the methods developed herein are only applicable to efficient evaluation of BIOs on prolate and/or oblate geometries. Owing to the linearity of the underlying PDEs, however, their utility clearly extends to particulate systems comprising of a mixture spheroidal and (possibly deformable) non-spheroidal particles. While the method developed here can be employed on the former, techniques such as pole-rotation-based quadrature \cite{veerapaneni2011fast} or MFS \cite{broms2024method} can be employed on the latter. 
 
In \cite{corona2018boundary}, a set of vector spherical harmonics were shown to diagonalize the Stokes BIOs for spherical particles. However, we have not yet found equivalent vector basis set for spheroidal particles; instead, we employ the standard approach of converting the Stokes BIOs to Laplace BIOs \cite{tornberg2008fmm3dstokes} and then apply scalar spheroidal harmonic formulas.

Lastly, we note that simulating large-scale dense particulate systems requires robust schemes for resolving collisions between particles that are inevitable due to discretization in space and time. In \cite{yan2019scalable}, a linear complementarity approach that enforces non-overlapping constraints at each time-step is employed to resolve collisions in a suspension of  spherical particles. Although the method is applicable to any convex particle shape, implementing the geometric constraints for spheroidal particles and incorporating it into the BIM developed here is an ongoing effort and will be reported at a later date. We refer the reader to the recent work of \cite{broms2024barrier} for a comprehensive review on this topic.   



\paragraph{Paper outline.} In Section~\ref{sec:background}, we introduce notation, spheroidal coordinate systems and bases required to represent integral densities on spheroidal surfaces and harmonic functions in the exterior and interior of spheroids. In Section~\ref{sec:potential}, we then derive evaluation formulas for Laplace layer potentials on and off one source spheroidal surface. We use these formulas in Section~\ref{sec:discretization} to develop a fast, spectrally accurate integral evaluation scheme for Laplace and Stokes BIOs. Finally, in Section~\ref{sec:numerical}, we perform a series of numerical experiments to validate our methods in the context of Laplace BVP solution and to test how accuracy and conditioning of commonly used BIOs depend on particle aspect ratio and proximity.

\section{Mathematical Preliminaries}
\label{sec:background}
In this section, we summarize the set of mathematical tools which we will use to represent scalar and vector fields on spheroidal surfaces. It is well-known that the Laplace operator is separable in spheroidal coordinates, and solutions to the Laplace equation can be written as an expansion in spheroidal harmonics bases on and off the spheroidal surface. As a continuation of the work in \cite{vico2014boundary,corona2018boundary,kohl2023fast}, we show that layer potential operators and their derivatives share this property, and so, we can derive analytical evaluation formulas on and off a spheroidal surface. 

\subsection{Spheroids and spheroidal coordinates}

We recall that a spheroid, or an ellipsoid of revolution, is a surface obtained by rotating an ellipse about one of its axes. As a result, it has a circular cross-section of radius $A$. In cartesian coordinates, points on an axis-aligned model surface satisfy the equation

\begin{equation*}
\frac{x_1^2+x_2^2}{A^2}+\frac{x_3^2}{C^2} = 1
\end{equation*}

If we rotate the ellipse about its major axis, the spheroid is prolate ($C>A$). Otherwise, $A$ is equal to the major semi-axis, and the spheroid is oblate ($C<A$). We will use the following definitions for prolate and oblate spheroidal coordinates:

\begin{center}
    \begin{tabular}{cc|c}
    &{\em \textbf{Prolate}} & {\em \textbf{Oblate}}  \\
    \hline 
    & &  \tabularnewline[-0.5em]
    \renewcommand{\arraystretch}{1.3}
    $\begin{pmatrix}
    x_1\\
    x_2\\
    x_3
    \end{pmatrix} = $ & 
    \renewcommand{\arraystretch}{1.3}
    $\begin{pmatrix}
    a\sqrt{u^2-1} \sqrt{1-v^2} \cos \phi\\
   a\sqrt{u^2-1} \sqrt{1-v^2} \sin \phi\\
    a u v
    \end{pmatrix}$ &
    \renewcommand{\arraystretch}{1.3}
    $\begin{pmatrix}
    a\sqrt{u^2+1} \sqrt{1-v^2} \cos \phi\\
    a\sqrt{u^2+1} \sqrt{1-v^2} \sin \phi\\
    a u v
    \end{pmatrix}$     \\    
    &  & \tabularnewline[-0.5em]
    \hline  
    &  & \tabularnewline[-0.5em]
    
    $\epsilon=$ & $1/u$ & $1/\sqrt{u^2+1}$ \\
    &  & \tabularnewline[-0.5em]
    & $u \in [1,\infty)$ &  $u \in [0,\infty)$ \\ 
    \end{tabular}
    \captionof{table}{{\em \small Coordinate and eccentricity of prolate and oblate spheroids.}}
 \label{tbl:coordinates}
\end{center} \vspace{0.05in}

with $\phi \in [0,2\pi]$ and $v \in [-1,1]$. Note that equivalent coordinate systems can be obtained by substituting $v = \cos \theta, u = \cosh \mu$ for prolates and $v = \sin \theta, u = \sinh \mu$ for oblates. In this coordinate system $(u,v,\phi)$, the surface defined by $u=u_0$ is a spheroid with semi-axes $A = a\sqrt{u_0^2 \pm 1}$ and $C = a u_0$. The magnitude of $u_0$ determines the eccentricity $\epsilon$ of the spheroid as indicated on Table \ref{tbl:coordinates}. Thus, when $u_0$ is close to its lower limit ($1$ for prolates, $0$ for oblates), the spheroid is more eccentric, but as $u_0 \to \infty$, the spheroid becomes more spherical. 

\subsection{Differential operators}

If we denote the correspondence with a point in cartesian coordinates $\mathbf{x} = \mathbf{x}(u,v,\phi)$, we can define the canonical coordinate vectors e.g. $\mathbf{e}_u = \partial_u \mathbf{x}$, $h_u = ||\mathbf{e}_u||$. Standard formulas for the gradient, divergence and other operators follow. We list here a few useful formulas for each spheroid type. 

\begin{center}
   \begin{tabular}{ c  c | c }
   & {\em \textbf{Prolate}} & {\em \textbf{Oblate}} \\
   \hline 
   & & \tabularnewline[-0.5em]
   \renewcommand{\arraystretch}{1.3}
     $\begin{pmatrix}
    h_\phi\\
    h_v \\
    h_u
    \end{pmatrix} = $ &
 \renewcommand{\arraystretch}{1.3}
     $\begin{pmatrix}
    a\sqrt{u^2-1}\sqrt{1-v^2} \\
    a\frac{\sqrt{u^2-v^2}}{\sqrt{1-v^2}} \\
    a\frac{\sqrt{u^2-v^2}}{\sqrt{u^2-1}}
    \end{pmatrix} $ 
     &  \renewcommand{\arraystretch}{1.3}
     $\begin{pmatrix}
   a\sqrt{u^2+1}\sqrt{1-v^2}\\
  a\frac{\sqrt{u^2+v^2}}{\sqrt{1-v^2}}\\
    a\frac{\sqrt{u^2+v^2}}{\sqrt{1+u^2}}
    \end{pmatrix}$ \\
    & &  \tabularnewline[-0.5em]
   \hline
   & & \tabularnewline[-0.5em]
    $\Delta f = $  & $\frac{\partial_u \left[ (u^2-1) f_u \right] + \partial_v \left[ (1-v^2) f_v \right]}{a^2(u^2-v^2)} + \frac{f_{\phi \phi}}{a^2(u^2-1)(1-v^2)} $ & $\frac{\partial_u \left[ (u^2+1) f_u \right] + \partial_v \left[ (1-v^2) f_v \right]}{a^2(u^2+v^2)}+ \frac{f_{\phi \phi}}{a^2(u^2+1)(1-v^2)}$  \\ 
   & & \tabularnewline[-0.5em]
   \hline
   & & \tabularnewline[-0.5em]
   $dS = $ & $a^2 \sqrt{u^2-v^2}\sqrt{u^2-1} \, dv\, d\phi$ & $a^2 \sqrt{u^2+v^2}\sqrt{u^2+1} \, dv\, d\phi$ \\
   
   \end{tabular}
   \captionof{table}{{\em \small Useful differential identities in the prolate and oblate spheroidal coordinates.}}
   \label{tab:differentials}
\end{center} \vspace{0.05in} 

\subsection{Spheroidal harmonic bases}

The spheroidal surface harmonic functions are, except for a change of variable (e.g. $v = \cos \theta$ for the prolate case) identical to their spherical counterparts. That is, 

\begin{definition} The scalar spheroidal harmonic $Y_n^m$ of degree $n$ and order $m$ (for $|m|\leq n$) is defined in terms of the associated Legendre functions $P_n^m$ by

\begin{equation*}
Y_n^m(v,\phi) = \sqrt{\frac{2n+1}{4\pi}} \sqrt{\frac{(n - m)!}{(n+m)!}} P_n^m(v) e^{i m \phi}
\end{equation*}

with $v \in [-1,1], \phi \in [0,2\pi]$. 
\end{definition}

For a given spheroid surface $\mathbb{S}_u$, spheroidal harmonics form an orthonormal basis of eigenfunctions of the Laplacian for all functions $\sigma \in L_w^2(\mathbb{S}_u)$, i.e. square integrable functions with an associated weighted inner product 

\begin{equation*}
\langle f,g \rangle_w = \int_{\mathbb{S}_u} f(v,\phi) g(v,\phi) w_u(v) dS  
\end{equation*}

with $w_u(v) = 1/ \sqrt{u^2-v^2}$ for prolates, and $w_u(v) = 1/\sqrt{u^2+v^2}$ for oblates. \\

For $\sigma \in C^{\infty}(\mathbb{S}_u)$, a classical result shows that truncating the expansion in spheroidal harmonics up to order $n = p$ yields a spectrally convergent approximation employing $(p+1)^2$ terms. The discrete spheroidal transform and its inverse then allow us to go to and from observations $\sigma(u,v)$ on a regular spheroidal grid, and a corresponding set of coefficients $\widehat{\sigma}_n^m$, respectively. While it is possible to improve the scaling of these from $O(p^4)$ to $O(p^2 \log^2 p)$ using Fast Fourier and Fast Legendre Transforms, break-even points for FLTs are typically large. For this reason, we employ only FFT accelerated spheroidal transforms with a complexity of $O(p^3 \log p)$. 

\subsubsection*{Solid spheroidal harmonics}

In order to evaluate solutions of the Laplace equation off spheroidal surfaces, we must find expansions in solid harmonics via separation of variables. For each spheroid type, we have a set of \emph{regular} solid spheroidal harmonics $R_n^m(u,v,\phi)$ and a set of \emph{irregular} solid spheroidal harmonics $I_n^m(u,v,\phi)$, constituting the solution space for the Laplace boundary value problem in the interior and exterior of the spheroid, respectively. Given a spheroid surface, corresponding to setting $u = u_0$, we evaluate solutions for the interior and exterior problems as an expansion of

\begin{center}
\begin{tabular}{ c  c | c | c }
& & {\em \textbf{Prolate}} & {\em \textbf{Oblate}} \\   
\hline
& & & \tabularnewline[-0.5em]
Interior & $R_n^m(u,v,\phi)$ & $P_n^m(u) Y_n^m(v,\phi)$ & $P_n^m(i u) Y_n^m(v,\phi)$ \\
& & & \tabularnewline[-0.5em]
Exterior & $I_n^m(u,v,\phi)$ & $Q_n^m(u) Y_n^m(v,\phi)$ & $Q_n^m(i u) Y_n^m(v,\phi)$ 
\end{tabular} \vspace{0.05in}
\end{center}

where $P_n^m , Q_n^m$ are associated Legendre functions of the first and second kind, respectively. 

\subsubsection*{Legendre function evaluation} 

The main computational task required to compute solid harmonic expansions is the evaluation of Legendre functions ``off-the-cut": for the prolate case, this involves $u > 1$, and for the oblate case, it involves $u>0$, as the argument $i u$ is on the positive imaginary semi-axis. For this purpose, we implement the algorithms in \cite{segura1999evaluation}, shown to be stable for prolate, oblate and toroidal harmonics. For each of these bases, there are two families $P_n^m,Q_n^m$ of solutions satisfying a forward recurrence relation. This recursion is a stable method to compute dominant solutions $P_n^m$, but is well-known to be unstable when applied to $Q_n^m$. Instead, the Wronskian relation is used to obtain a stable backwards recursion as follows: 

\begin{enumerate}
   \item The forward recurrence relation and formulas for $n=m,m+1$:
    \begin{align}
    &(n-m+1)P_{n+1}^m(x)- (2n+1)xP_n^m(x)+(n+m)P_{n-1}^m = 0, \quad n\geq 1 \label{eq:recur} \\ 
    &P_m^m(x)=(2m-1)!! (x^2-1)^{m/2} \ \ \ 
    P_{m+1}^m(x)=x(2m+1)P_m^m(x)
    \end{align}
   are used to compute $P_n^m$ for $n$ and $|m|\leq n$ up to a desired order $N$. 
   \item A continued fraction formula is derived from Eq. \eqref{eq:recur}: 

   \begin{equation}\label{eq:continuedfraction}
    H_n^m(x)=\frac{Q_n^m(x)}{Q_{n-1}^m(x)}=\frac{-\frac{n+m}{n-m+1}}{-\frac{2n+1}{n-m+1}x + \frac{Q_{n+1}^m(x)}{Q_n^m(x)}} = \frac{-\frac{n+m}{n-m+1}}{-\frac{2n+1}{n-m+1}x + H_{n+1}^m(x)}, \quad n\geq 1
\end{equation}

A modified Lentz algorithm \cite{press} is used to compute $H_N^m(x)$ using Eq \eqref{eq:continuedfraction}.

   \item We combine Eq \eqref{eq:continuedfraction} with the Wronskian relation

   \begin{equation}\label{eq:PQ}
    P_n^m(x) Q_{n-1}^m(x) - P_{n-1}^m(x) Q_n^m(x) =  \frac{(n+m-1)!}{(n-m)!}(-1)^m, \quad n\geq 1
    \end{equation}
    to obtain 
    \begin{align}
    Q_{N-1}^m(x) &= \frac{(n+m-1)!}{(n-m)!} \frac{(-1)^m}{P_N^m(x)-H_N^m(x)P_{N-1}^m(x)}\label{eq:QfromP1}\\
    Q_N^m(x) &= H_N^m(x)Q_{N-1}^m(x)\label{eq:QfromP2}
    \end{align}
    \item We then use Eqs \eqref{eq:QfromP1},\eqref{eq:QfromP2} to compute $Q_N^m,Q_{N-1}^m$, and recur backwards to $Q_m^m$.
\end{enumerate}

\vskip 5 pt

Accurate and stable evaluation of first derivatives $P_n^{m \prime}$ and $Q_n^{m \prime}$ is also required in the layer potential evaluation formulas in Section \ref{sec:potential}. For this purpose, it suffices to use 

\begin{equation}\label{eq:derivative}
    (1-x^2) P_n^{m \prime}(x)=(m-n-1)P_{n+1}^m(x) +(n+1)x P_n^m(x)
\end{equation}

which is satisfied by both families of Legendre functions. 

\subsection{Boundary integral operators}\label{ssc:bio}

Layer potential operators are used to represent solutions of the Laplace equation in terms of a combination of unknown charge and dipole densities defined on domain boundaries; this then leads to boundary integral equations for these densities. Let $\Gamma=\PS{u_0}$ denote the surface of a prolate spheroid defined by $u=u_0>1$ in spheroidal coordinates. We denote the free space fundamental solution for the Laplace equation $G(\bm{x},\bm{y}) = \frac{1}{4 \pi ||\bm{x}-\bm{y}||}$, and $\nu(\bm{y})$ the outward normal vector at $\bm{y} \in \PS{u_0}$. The Laplace single and double layer potential operators are then 
\begin{align}
\SLP{u_0}[\sigma](\bm{x}) &= \int_{\PS{u_0}} G(\bm{x},\bm{y}) \sigma(\bm{y}) dS(\bm{y}) \label{eq:SLP} \\
\DLP{u_0}[\mu](\bm{x}) &= \int_{\PS{u_0}} \frac{\partial G}{\partial \bm{\nu}(\bm{y})}(\bm{x},\bm{y}) \mu(\bm{y}) dS(\bm{y}) \label{eq:DLP}
\end{align}
Similarly, we will use the notation $\SLO{u_0}[\sigma], \DLO{u_0}[\mu]$ for single and double layer operators defined on $\Gamma=\OS{u_0}$, the surface of an oblate spheroid defined by $u_0>0$ in oblate spheroidal coordinates. We assume the integral densities $\sigma,\mu$ are in $L^2(\Gamma)$. By construction, these layer potentials are solutions to the Laplace equation both in the exterior and interior of the spheroid; we can then take the appropriate limit as $x \rightarrow \Gamma$ along the normal direction to define operator evaluation on the surface. 

\subsubsection*{Jump conditions}

We briefly summarize standard jump conditions \cite{kress} for Laplace layer potentials. That is, if $\varphi$ is a layer potential or its normal derivative applied to a density (e.g. $\sigma, \mu$), we are interested in the difference between exterior and interior limits, denoted as $[ [ \varphi ] ]_\Gamma = \varphi^e - \varphi^i$. 

\begin{align*}
[[ \SL[\sigma] ]]_\Gamma &= 0 \ \ \ \left[ \left[ \frac{\partial}{\partial \bm{\nu}} \SL[\sigma] \right] \right]_\Gamma = -\sigma \\
[[ \DL[\mu] ]]_\Gamma &= \mu \ \ \ \left[ \left[ \frac{\partial}{\partial \bm{\nu}} \DL[\mu] \right] \right]_\Gamma = 0
\end{align*}

For those potentials like the single layer which are continuous across $\Gamma$, we use the same notation for the operator on the exterior, interior and on-surface (singular). For the double layer and the normal derivative of the single layer, we have that: 
\begin{align*}
\lim_{x \rightarrow \Gamma^\pm} \DL_{u_0}[\mu](\bm{x}) &= \DL_{u_0}^\pm[\mu](\bm{x}) = \DL_{u_0}[\mu](\bm{x}) \pm \frac{1}{2} \mu(\bm{x}), \\
\lim_{x \rightarrow \Gamma^\pm} \frac{\partial}{\partial \bm{\nu}(\bm{x})}\SL_{u_0}[\sigma](\bm{x}) &= (\SL')^\pm_{u_0}[\sigma](\bm{x}) = \SL'_{u_0}[\sigma](\bm{x}) \mp \frac{1}{2} \sigma(\bm{x}).
\end{align*}
here $\DL_{u_0}[\mu], \SL'_{u_0}[\sigma]$ denote singular integral operators, e.g. Eq. \ref{eq:DLP} for $\bm{x}\in \Gamma$.

\section{Layer potential evaluation}
\label{sec:potential}

Given an integral density $\sigma \in L^2(\Gamma)$, it can be expanded in the spheroidal harmonics basis as
\begin{equation}
\sigma(v,\phi) = \sum_{n=0}^\infty \sum_{m=-n}^n \widehat{\sigma_n^m} Y_n^m(v,\phi) \label{eq:harmexp}
\end{equation}
We then propose a representation of the solution to a Laplace boundary value problem in terms of a linear combination of layer potentials and their derivatives. Say we select $\varphi(x) = \SL[\sigma](x)$. This function is, by construction, a solution to the Laplace equation for all $x$ except on the surface $\Gamma$. Since $\varphi$ is a solution to the Laplace equation in the interior or the exterior of this spheroid, it  can be written as an infinite expansion in regular or irregular solid harmonics, of the form: 

\begin{equation}
\varphi(u,v,\phi) = \sum_{n=0}^{\infty} \sum_{m=-n}^n a_n^m f^m_n(u) Y_n^m(v,\phi) \label{eq:sharmexp}
\end{equation}

For both the interior and exterior cases, we now have a way to evaluate the one-sided limits as $x \rightarrow \Gamma^{\pm}$ in the normal direction: by taking the limit as $u \rightarrow u_0^{\pm}$ in Eq. \eqref{eq:sharmexp}. We can then use standard jump conditions for $\varphi$ and $\frac{\partial \varphi}{\partial \nu}$ across $\Gamma$. Finally, using orthogonality of the basis $Y_n^m$, this results in a system of two linear equations on the unknowns $a_n^{m,E}$ and $a_n^{m,I}$ for each $(n,m)$ pair, which can be solved analytically.   

\subsection{Double layer potential evaluation}\label{sec:DL_eval}

We begin by applying the derivation procedure described above to the double layer because it readily diagonalizes in the spheroidal harmonic basis, both on and off the surface. In other words, given $\mu \in L^2(\Gamma)$, we can find expansions for $\varphi = \DL[\mu]$ on (Eq. \eqref{eq:harmexp}) and off (Eq. \eqref{eq:sharmexp}) the surface $\Gamma$ where the $(n,m)$ coefficient is $\lambda_n^m \widehat{\mu_n^m}$. 

\begin{theorem}[Spheroidal Double-layer potential] Let $\Gamma$ be the surface of a prolate spheroid $\PS{u_0}$ ($u_0>1$) or an oblate spheroid $\OS{u_0}$ ($u_0>0$), and let $\mu \in L^2(\Gamma)$ with spheroidal harmonic coefficients $\widehat{\mu_n^m}$. Then, we have the following evaluation formulas: 

\begin{center}
   \begin{tabular}{c|c|c}
       & \textbf{Prolate} $\DLP{u_0}[\mu]$ & \textbf{Oblate} $\DLO{u_0}[\mu]$ \\
       \hline
       & & \tabularnewline[-0.5em]
       Exterior & $\sum\limits_{n,m} b_n^m \widehat{\mu_n^m} P_n^{m \prime}(u_0)Q_n^m(u) Y_n^m(v,\phi)$ &  
       $\sum\limits_{n,m} c_n^m \widehat{\mu_n^m} P_n^{m \prime}(i u_0)Q_n^m(i u) Y_n^m(v,\phi)$ \\
       & & \tabularnewline[-0.5em]
       Interior & $\sum\limits_{n,m} b_n^m \widehat{\mu_n^m} Q_n^{m \prime}(u_0) P_n^m(u) Y_n^m(v,\phi)$ & 
       $\sum\limits_{n,m} c_n^m \widehat{\mu_n^m} Q_n^{m \prime}(i u_0) P_n^m(i u) Y_n^m(v,\phi)$
   \end{tabular}
\end{center}\vspace{0.05in}

where exterior and interior correspond to $u>u_0$ and $u<u_0$, respectively. $b_n^m, c_n^m$ are given by: 

\begin{equation*}
b_n^m = \frac{(n-m)!}{(n+m)!}(-1)^m(u_0^2-1) \qquad \qquad c_n^m = \frac{(n-m)!}{(n+m)!}(-1)^{m+1}(u_0^2+1) 
\end{equation*}

Spectra for $\DL^+$ and $\DL^-$ can be obtained by plugging $u=u_0$: 

\begin{center}
   \begin{tabular}{c|c|c}
       & \textbf{Prolate} $\DLP{u_0}[\mu]$ & \textbf{Oblate} $\DLO{u_0}[\mu]$ \\
       \hline
       & & \tabularnewline[-0.5em]
       $\DL^+$ & $\sum\limits_{n,m}b_n^m \widehat{\mu_n^m} P_n^{m \prime}(u_0)Q_n^m(u_0)$ &  
       $\sum\limits_{n,m}c_n^m \widehat{\mu_n^m} P_n^{m \prime}(i u_0)Q_n^m(i u_0)$ \\
       & & \tabularnewline[-0.5em]
       $\DL^-$ & $\sum\limits_{n,m}b_n^m \widehat{\mu_n^m} Q_n^{m \prime}(u_0) P_n^m(u_0)$ & 
       $\sum\limits_{n,m}c_n^m \widehat{\mu_n^m} Q_n^{m \prime}(i u_0) P_n^m(i u_0)$
   \end{tabular}
\end{center} \vspace{0.05in}

The spectra for the singular integral operator $\DL$ is then the average of that for $\DL^+$ and $\DL^-$. 
\label{thm:DLthm}
\end{theorem}

\vskip 5 pt

We include here the derivation for the prolate case, and then indicate what changes are needed for the oblate case. Applying the jump conditions in Section \ref{ssc:bio} to the double-layer given by
\begin{align*}
    \DLP{u_0}[\mu]^-(u,v,\phi)= \sum_{m,n} a_n^{m,I} P_n^m(u)Y_n^m(v,\phi) \quad & \quad \text{(interior)}\\
    \DLP{u_0}[\mu]^+(u,v,\phi)=\sum_{m,n} a_n^{m,E} Q_n^m(u)Y_n^m(v,\phi) \quad & \quad \text{(exterior)}
\end{align*}
gives us equation \eqref{eq:DLjumpeq}. We then apply the normal derivative operator in prolate spheroidal coordinates given by $\frac{\partial}{\partial \nu}=\frac{1}{a}\sqrt{\frac{u^2-1}{u^2-v^2}}\frac{\partial}{\partial u}$. Using continuity of $\DL'$ across $\Gamma$, we obtain \eqref{eq:DLprimejumpeq}. 
\begin{align}
    &\sum_{m,n} a_n^{m,E} Q_n^m(u_0)Y_n^m(v,\phi)-\sum_{m,n} a_n^{m,I} P_n^m(u_0)Y_n^m(v,\phi) = \sum_{n,m} \widehat{\mu_n^m} Y_n^m(v,\phi) \label{eq:DLjumpeq} 
    \\
    &\sum_{m,n} a_n^{m,E} Q_n^{m\prime}(u_0)Y_n^m(v,\phi)-\sum_{m,n} a_n^{m,I} P_n^{m\prime}(u_0)Y_n^m(v,\phi)=0 \label{eq:DLprimejumpeq}
\end{align}
%
%
%
Using orthogonality of spheroidal harmonics and Eqs \eqref{eq:DLjumpeq} 
and \eqref{eq:DLprimejumpeq}
, we obtain the system: 
\begin{equation}
    \bma Q_n^m(u_0) & -P_n^m(u_0)\\ Q_n^{m\prime}(u_0) & -P_n^{m\prime}(u_0) \ema \bma a_n^{m,E} \\ a_n^{m,I} \ema = \bma \widehat{\mu_n^m} \\ 0\ema
\end{equation}

The solution of this system yields the evaluation formulas for the prolate case; we note that the determinant of this matrix is the Wronskian $W(P_n^m,Q_n^m)(u_0) = \frac{(n+m)!(-1)^m}{(n-m)!(u_0^2-1)}$. The derivation of the oblate case is practically identical, except that $\frac{\partial}{\partial \nu}=\frac{1}{a}\sqrt{\frac{u^2+1}{u^2+v^2}}\frac{\partial}{\partial u}$, inputs on all Legendre functions are of the form $i u$, and the formula for the Wronskian changes to $W(P_n^m,Q_n^m)(i u_0) = \frac{(n+m)!(-1)^{m+1}}{(n-m)!(u_0^2+1)}$.

\subsection{Single layer potential evaluation}

Unlike the double layer, the single layer potential does not diagonalize on spheroidal coordinates. Following the steps for the derivation for the double layer, it becomes apparent that this is because the normal derivative operator depends on $v$ (and hence the azimuthal angle). This introduces a factor of the form $1/\sqrt{u_0^2 \pm v^2 }$ on the jump conditions for $\SL'$. 

However, this in turn tells us that the modified potentials $\widetilde{\SLP{u_0}}[\sigma] = \SLP{u_0}[\sigma / \sqrt{u_0^2-v^2}]$ and $\widetilde{\SLO{u_0}}[\sigma] = \SLO{u_0}[\sigma / \sqrt{u_0^2+v^2}]$ diagonalize with respect to the respective spheroidal harmonic bases, and we can apply the same derivation procedure as above. Given $\sigma \in L^2(\Gamma)$, we denote as $\widetilde{\sigma_n^m}$ the coefficients of the expansion of $\sigma(v,\phi)$ on the basis $\left\{ \frac{Y_n^m}{\sqrt{u_0^2 \pm v^2}} \right\}$. That is, 

\begin{equation}
\sigma(v,\phi) = \sum_{n=0}^{\infty} \sum_{m=-n}^n \widetilde{\sigma_n^m} \frac{Y_n^m(v,\phi)}{\sqrt{u_0^2 \pm v^2}} 
\end{equation}

We note that since $v \in [-1,1]$, the weight $(u_0^2 \pm v^2)^{-1/2}$ for both prolate and oblate cases is always positive and in $C^{\infty}(\Gamma)$. This means we can safely factor it in and out of expansions like the one above, and find the coefficients $\widetilde{\sigma_n^m}$ via a truncated expansion in spherical harmonics of $(u_0^2 \pm v^2)^{1/2} \sigma(v,\phi)$.

\vskip 5 pt 

\begin{theorem}[Spheroidal Single-layer potential]
Let $\Gamma$ be the surface of a prolate spheroid $\PS{u_0}$ ($u_0>1$) or an oblate spheroid $\OS{u_0}$ ($u_0>0$), and let $\sigma \in L^2(\Gamma)$ with coefficients $\widetilde{\sigma_n^m}$ in the modified basis. Then, we have the following evaluation formulas: 

\begin{center}
   \begin{tabular}{c|c|c}
       & \textbf{Prolate} $\SLP{u_0}[\sigma]$ & \textbf{Oblate} $\SLO{u_0}[\sigma]$ \\
       \hline
       & & \tabularnewline[-0.5em]
       Exterior & $\sum\limits_{n,m} \widetilde{b_n^m} \widetilde{\sigma_n^m} P_n^{m}(u_0)Q_n^m(u) Y_n^m(v,\phi)$ &  
       $\sum\limits_{n,m}  \widetilde{c_n^m} \widetilde{\sigma_n^m} P_n^{m}(i u_0)Q_n^m(i u) Y_n^m(v,\phi)$ \\
       & & \tabularnewline[-0.5em]
       Interior & $\sum\limits_{n,m} \widetilde{b_n^m} \widetilde{\sigma_n^m} Q_n^{m}(u_0) P_n^m(u) Y_n^m(v,\phi)$ & 
       $\sum\limits_{n,m} \widetilde{c_n^m} \widetilde{\sigma_n^m} Q_n^{m}(i u_0) P_n^m(i u) Y_n^m(v,\phi)$
   \end{tabular}
\end{center}\vspace{0.05in}

where exterior and interior correspond to $u>u_0$ and $u<u_0$, respectively. 
$\widetilde{b_n^m}, \widetilde{c_n^m}$ are given by: 

\begin{equation*}
\widetilde{b_n^m} = a\frac{(n-m)!}{(n+m)!}(-1)^m \sqrt{u_0^2-1} \qquad \qquad 
\widetilde{c_n^m} = ia \frac{(n-m)!}{(n+m)!}(-1)^{m}\sqrt{u_0^2+1} 
\end{equation*}

Surface spheroidal coefficients result from evaluating $u=u_0$ ($\SL$ is continuous, so $\SL^+ = \SL^- = \SL$):

\begin{center}
   \begin{tabular}{c|c|c}
       & \textbf{Prolate} $\SLP{u_0}[\sigma]$ & \textbf{Oblate} $\SLO{u_0}[\sigma]$ \\
       \hline
       & & \tabularnewline[-0.5em]
       $\SL$ & $\sum\limits_{n,m}\widetilde{b_n^m} \widetilde{\sigma_n^m} P_n^{m}(u_0)Q_n^m(u_0)$ &  
       $\sum\limits_{n,m}\widetilde{c_n^m} \widetilde{\sigma_n^m} P_n^{m}(i u_0)Q_n^m(i u_0)$ 
   \end{tabular}
\end{center} \vspace{0.05in}

\label{thm:SLthm}
\end{theorem} 

\vskip 5 pt

Once again, it suffices to include the derivation for the prolate case. We apply jump conditions to

\begin{align*}
    \SLP{u_0}[\sigma]^-(u,v,\phi)= \sum_{m,n} \widetilde{a_n^{m,I}} P_n^m(u) Y_n^m(v,\phi) \quad & \quad \text{(interior)}\\
    \SLP{u_0}[\sigma]^+(u,v,\phi)=\sum_{m,n} \widetilde{a_n^{m,E}} Q_n^m(u) Y_n^m(v,\phi) \quad & \quad \text{(exterior)}
\end{align*}

where now the single layer is continuous across $\Gamma$, and its normal derivative has a non-trivial jump. We arrive at the following equations: 

\begin{align}
&\sum_{m,n} \widetilde{a_n^{m,E}} Q_n^m(u_0)Y_n^m(v,\phi)-\sum_{m,n} \widetilde{a_n^{m,I}} P_n^m(u_0)Y_n^m(v,\phi) = 0 \label{eq:SLjumpeq} \\
    &\frac{1}{a}\sqrt{u_0^2-1} \left[\sum_{m,n} \widetilde{a_n^{m,E}} Q_n^{m\prime}(u_0)Y_n^m(v,\phi)-\sum_{m,n} \widetilde{a_n^{m,I}} P_n^{m\prime}(u_0)Y_n^m(v,\phi)\right]= -\sum_{n=0}^{\infty} \sum_{m=-n}^n \widetilde{\sigma_n^m} Y_n^m(v,\phi) \label{eq:SLprimejumpeq}
\end{align}

Using orthogonality of spheroidal harmonics and Eqs \eqref{eq:SLjumpeq} and \eqref{eq:SLprimejumpeq}, we obtain the system: 

\begin{equation}
    \bma Q_n^m(u_0) & -P_n^m(u_0)\\ Q_n^{m\prime}(u_0) & -P_n^{m\prime}(u_0) \ema \bma \widetilde{a_n^{m,E}} \\ \widetilde{a_n^{m,I}} \ema = \bma 0 \\ -\frac{a \widetilde{\sigma_n^m}}{\sqrt{u_0^2-1}} \ema
\end{equation}

formulas for these coefficients follow from the solution of this $2 \times 2$ system. We note the system matrix is the same as for the double layer. 

\subsection{Derivatives of layer potentials}\label{sec:derivatives}
The normal derivative operator of the single layer potential, denoted $\SL^{\prime}$, is useful when solving Neumann boundary value problems, and is the adjoint operator of the double layer potential \cite{folland}. We derive evaluation formulas on and off-surface from single layer potential evaluations in the prolate case, and indicate changes to be made in the case of oblates. We note that the evaluation formula and procedure below can be straight-forwardly generalized to any directional derivative of layer potentials, e.g. $\DL^{\prime}$. 

\begin{theorem}[Derivative of the Spheroidal Single-layer potential]
Let $\Gamma$ be the surface of a prolate spheroid $\PS{u_0}$ ($u_0>1$) or an oblate spheroid $\OS{u_0}$ ($u_0>0$), and let $\sigma \in L^2(\Gamma)$ with coefficients $\widetilde{\sigma_n^m}$ in the modified basis. Then, we have the following evaluation formulas on $\Gamma$: 

\begin{center}
   \begin{tabular}{c|c|c}
       & \textbf{Prolate} $\SLPp{u_0}[\sigma]$ & \textbf{Oblate} $\SLOp{u_0}[\sigma]$ \\
       \hline
       & & \tabularnewline[-0.5em]
       $\mathcal{S^{\prime+}}$ & $\sum\limits_{n,m}\frac{\widetilde{b_n^m}\sqrt{u_0^2-1}}{a} \widetilde{\sigma_n^m} P_n^{m}(u_0)Q_n^{m\prime}(u_0) \frac{Y_n^m(v,\phi)}{\sqrt{u_0^2-v^2}}$  & $\sum\limits_{n,m}i\frac{\widetilde{c_n^m}\sqrt{u_0^2+1}}{a} \widetilde{\sigma_n^m} P_n^{m}(i u_0)Q_n^{m\prime}(i u_0)\frac{Y_n^m(v,\phi)}{\sqrt{u_0^2+v^2}}$ \\
       $\mathcal{S^{\prime-}}$ & $\sum\limits_{n,m}\frac{\widetilde{b_n^m}\sqrt{u_0^2-1}}{a} \widetilde{\sigma_n^m} Q_n^{m}(u_0)P_n^{m\prime}(u_0)\frac{Y_n^m(v,\phi)}{\sqrt{u_0^2-v^2}}$ & $\sum\limits_{n,m}i\frac{\widetilde{c_n^m}\sqrt{u_0^2+1}}{a} \widetilde{\sigma_n^m} Q_n^{m}(i u_0)P_n^{m\prime}(i u_0)\frac{Y_n^m(v,\phi)}{\sqrt{u_0^2+v^2}}$
   \end{tabular}
\end{center}\vspace{0.05in}
Coefficients for the derivative of operator $\SL^{\prime}$ on surface are the average of those for $\SL^{\prime+}$ and $\SL^{\prime-}$.
\label{thm:dSLthm}
\end{theorem}

Given $\sigma\in L^2(\Gamma)$, Theorem \ref{thm:SLthm} tells us $S^{\mathbb{P}}_{u_0}[\sigma](u,v,\phi)=\sum\limits_{n,m}\widetilde{b_n^m} \widetilde{\sigma_n^m} g_n^mf_n^m(u)Y_n^m(v,\phi)$, where for prolates, $g_n^m=P_n^m(u_0)$, $f_n^m=Q_n^m$ on the exterior and  $g_n^m=Q_n^m(u_0)$, $f_n^m=P_n^m$ on the interior. For ease of computation, we denote the normalization factor in Definition 1 as $N_n^m$, and so,  
$Y_n^m(v,\phi)=N_n^mP_n^m(v)e^{im\phi}$. Given the unit normal vector $\bm\nu$ at a target point $x$, with $\bm\nu=\nu_u\bm e_u+\nu_v \bm e_v+\nu_{\phi}\bm e_{\phi}$,
\begin{align}
    \frac{\partial}{\partial\bm{\nu}}\mathcal{S}\left[\sigma\right] &= \nabla \left[\sum_{n,m}\widetilde{b_n^m} \widetilde{\sigma_n^m} g_n^m f_n^m(u)Y_n^m(v,\phi)\right]\cdot \bm \nu \notag \\ 
    &= \sum_{n,m} \widetilde{b_n^m} \widetilde{\sigma_n^m} g_n^m \Bigg[\frac{1}{a} \sqrt{\frac{u^2-1}{u^2-v^2}} \nu_u  f_n^{m\prime}(u)N_n^mP_n^m(v)e^{im\phi} \notag \\
    &\quad +  \frac{1}{a}\sqrt{\frac{1-v^2}{u^2-v^2}}\nu_v f_n^m(u) N_n^mP_n^{m\prime}(v)e^{im\phi} \notag \\
    &\quad + \frac{1}{a}\frac{1}{\sqrt{(u^2-1)(1-v^2)}}\nu_\phi  f_n^m(u)N_n^mP_n^m(v)(im)e^{im\phi}\Bigg]
\end{align}
Using the recurrence relation for derivatives of Legendre functions in Eq \eqref{eq:derivative}, we obtain 

\begin{align}
    \frac{\partial}{\partial\bm{\nu}}\mathcal{S}\left[\sigma\right] &= \sum_{n,m} \frac{\widetilde{b_n^m}}{a} \widetilde{\sigma_n^m} g_n^m \Bigg[ \left(\sqrt{\frac{u^2-1}{u^2-v^2}}\nu_uf_n^{m\prime}(u)+ \left(\frac{(n+1)v\nu_v }{\sqrt{u^2-v^2}}+ \frac{im\nu_\phi}{\sqrt{u^2-1}}\right)\frac{f_n^m(u)}{\sqrt{1-v^2}}\right) Y_n^m(v,\phi) \notag \\
    -& \sqrt{\frac{2n+1}{2n+3}\frac{n+m+1}{n-m+1}}\frac{n-m+1}{\sqrt{(u^2-v^2)(1-v^2)}}\nu_v f_n^m(u)Y_{n+1}^m(v,\phi)\Bigg]
    \label{eq:dSLevalformula}
\end{align}

For a target $x$ approaching $\Gamma$ in the normal direction, we must take the limit $u \rightarrow u_0$, and take $\nu(x) = e_u(x)$ (that is, $\nu_u=1, \nu_v=\nu_\phi=0$). This yields expansions in spheroidal bases for $\mathcal{S^{\prime+}}$, $\mathcal{S^{\prime -}}$ 
taking the appropriate one-sided limit. However, Equation \eqref{eq:dSLevalformula} can be used to take arbitrary directional derivatives of the single layer operator. This is essential for near-field neighbor interactions and implementing the Stokes operators using Laplace operators.

For the oblate case, $S^{\mathbb{O}}_{u_0}[\sigma](u,v,\phi)=\widetilde{c_n^m}g_n^mf_n^m(iu)Y_n^m(v,\phi)$, where $g_n^m=P_n^m(iu_0)$, $f_n^m=Q_n^m$ on the exterior and  $g_n^m=Q_n^m(iu_0)$, $f_n^m=P_n^m$ on the interior. The gradient is also slightly different, as indicated in Table \ref{tab:differentials}. 

\section{Numerical scheme}
\label{sec:discretization}
In Section \ref{sec:potential}, we have derived evaluation formulas for the Laplace layer potentials $\DL$, $\SL$ and derivatives like $\SL'$, and $\DL'$ on and off the surface of the spheroid in terms of expansions in spheroidal harmonics bases.  In this section, we develop a fast, spectrally accurate evaluation scheme for layer potentials defined on a collection of spheroidal surfaces. We employ the formulas mentioned above to perform singular and near-singular integral evaluation. 

\subsubsection*{Problem setup} Given $M$ spheroidal particles (prolates and oblates) with surfaces $\Gamma_i$, we wish to evaluate 
\begin{equation}
u(\bm{x}) = \K_{\Gamma}[\sigma] (\bm{x})  = \sum_{i=1}^M \K_{\Gamma_i}[\sigma_i] (\bm{x})
\label{eq:potall}
\end{equation}
where $\sigma_i \in L^2(\Gamma_i)$ and each $\K_{\Gamma_i}$ is an integral operator of the form
\begin{equation}
    u_i(\bm{x}) = \K_{\Gamma_i}[\sigma_i] (\bm{x}) = \int_{\Gamma_i}F_i( \bm{x}, \bm{y}) \sigma_i( \bm{y})\rmd S(\bm{y}) \label{eq:poti}
\end{equation} 
with the integral kernel $F_i$ a linear combination of layer potentials or of their derivatives. In developing our evaluation scheme, we consider two important cases: 

\begin{enumerate}
\item \emph{Evaluation of $u(\bm{x})$ at all particle surfaces $\Gamma_i$.} This routine can then be used to solve boundary integral equations defined on $\Gamma = \cup_{i=1}^M \Gamma_i$, e.g. using Krylov subspace methods like GMRES. For instance, Algorithm \ref{alg:DL} performs this operation for the Laplace Double Layer potential.  
\item \emph{Evaluation of $u(\bm{x})$ at a given set of target points.} Given the integral densities $\sigma_i$ obtained from this BIE solver, we may also use the techniques described below to evaluate the solution to our PDE problem on any set of targets in the domain. 
\end{enumerate}
\subsubsection*{Integral density representation}

We consider two approximate representations for a scalar density $\sigma(\bm{x})$ defined on a spheroidal surface $\Gamma_i$. For a given order $p$, we approximate $u$ using a truncated expansion in spheroidal harmonics ($n \leq p, |m|\leq n$). This involves computing $(p+1)^2$ coefficients $\widehat{\sigma_n^m}$. Alternatively, we sample points $\m y_{j,k}=\m y(v_j,\phi_k$) on $\Gamma_i$, where $v_j$ are the $(p+1)$ Gauss-Legendre nodes in $[-1,1]$ and
$\phi_k$ are $2p$ equispaced discretization points in $[0,2\pi]$. This results in $2p(p+1)$ evaluations $\sigma(y_{j,k})$. We employ fast spherical harmonics transforms to go from one discrete representation to the other, as discussed in Section \ref{sec:background}. We note that our evaluation scheme for the single layer potential and its derivatives will require us to compute coefficients for modified bases, as indicated in Theorem \ref{thm:SLthm}. 

\subsection{Numerical integration regimes}

When we discretize our BIO, evaluation of $\K_{\Gamma_i}[\sigma_i] (\bm{x})$ involves target points that are either on the ``source'' surface $\Gamma_i$ or off surface. Because the integral kernel $F_i$ is singular when $\bm{x} \in \Gamma_i$, we must employ different techniques to accurately and efficiently evaluate these integrals in three different regimes: when targets are far from the the source (\textit{smooth}), when the targets are on the surface (\textit{singular}), and when they are off the surface, but close to it (\textit{near-singular}). To be precise, given a user-prescribed parameter $\eta \geq 0$, we will say $\bm{x}$ is ``far'' from $\Gamma_i$ (and part of its far-field) if 
\begin{equation}
    d(\bm{x},\Gamma_i) \geq \eta \ \mathrm{diam}(\Gamma_i);
\end{equation}
otherwise, $\bm{x}$ is ``near'' (in the near-field). We assume $\eta=1$ unless otherwise specified; this parameter can be tuned to meet target accuracy requirements. In our proposed method, off-surface targets are split into near and far fields, where far-field targets are evaluated with smooth quadrature, and near targets are evaluated using formulas based on spheroidal harmonics expansions.

\subsubsection*{Far-field: smooth integration}

We have chosen our grid on the $(v,\phi)$ spheroidal coordinates in accordance with a standard tensor quadrature scheme, employing Gauss-Legendre weights in $v$ and composite trapezoidal in $\phi$. This rule is known to be spectrally convergent for smooth integrands. Given a total of $N = O(Mp^2)$ discretization points on our spheroid suspension, the cost of directly evaluating all far-field interactions (a necessary computation for BIO matrix-vector apply as in  Algorithm \ref{alg:DL}) is $\mathcal{O}(N^2) = \mathcal{O}(M^2 p^4)$. Since this operation is a summation of Green’s functions or their derivatives, Fast Multipole Method (FMM) acceleration can be employed to reduce the cost of this operation to $\mathcal{O}(N) = \mathcal{O}(Mp^2)$.

\subsubsection*{Near-field: singular and near-singular integration}

Evaluating $\mathcal{K}_{\Gamma_i}$ at points belonging to the source spheroid $\Gamma_i$, that is, computing particle \textit{self-interactions}, requires evaluating an integral with a known singularity in the integrand. Thus, standard smooth quadrature methods such as those used for the far-field above will fail in this region. We also generally need to evaluate $\mathcal{K}_{\Gamma_i}$ at nearby off-surface points, such as would be present in interaction between the source particle and close neighbor particles. In this regime, a standard smooth quadrature is inefficient, requiring very fine discretization to accurately compute the near-singular integral. 

Evaluation of layer potentials and related BIOs in general 3D geometries requires the use of specialized techniques for singular and near-singular integration to accurately perform near-field evaluation. One of the key features of our scheme is that it relies on the analytical formulas derived in Section \ref{sec:potential} to 
evaluate both on-surface and near off-surface interactions via surface and solid spheroidal harmonic expansions, respectively.

\subsubsection*{Computational complexity of near-field interactions}

Given arbitrarily chosen $M_{trg}$ target points, 
we must for each point $\bm{x}$ evaluate and add $\mathcal{K}_{\Gamma_i}$ for all $i$ such that $\bm{x}$ is in $\Gamma_i$'s near-field. In most practical scenarios, we can bound the number of particles near one point by a constant $M_{neigh}$. Direct evaluation of the $O(p^2)$ terms in the corresponding spheroidal expansion formulas is then $\mathcal{O}(M_{trg} M_{neigh} p^2)$.

For all-to-all particle interactions as those involved in Algorithm \ref{alg:DL}, however, we must consider the $\mathcal{O}(Mp^2)$ target points corresponding to our discretization (or equivalently, $\mathcal{O}(Mp^2)$ spheroidal harmonics coefficients for $u(\bm{x})$). Following the reasoning above, direct evaluation of near-field interactions is then $\mathcal{O}(M M_{neigh} p^4) = \mathcal{O}(M_{neigh} p^2 N)$, where $M_{neigh}$ is now the maximum number of neighbors a given particle can have. In \cite{corona2018boundary}, accelerations are proposed for the case of spherical suspensions; as we attempt to extend these to the spheroidal case, we note that: 
\begin{itemize}
\item \emph{Self-interactions.} For layer potentials and their derivatives, the linear operator mapping coefficients $\widehat{\sigma_n^m}$ to $\widehat{u_n^m}$ (in the appropriate bases) is \emph{diagonal}. Evaluation at all target points on the spheroidal particle then involves a fast, $\mathcal{O}(p^3 \log p)$ transform. Overall, self-interactions can be computed in $\mathcal{O}(Mp^3 \log p) = \mathcal{O}(N p \log p)$ work.   
\item \emph{Neighbor particle interactions.} In this case, the ``point-and-shoot" FFT-based algorithm in \cite{corona2018boundary} cannot be directly extended, as spheroidal surface grids generally can not be aligned due to particle orientations. Implementation of general spheroidal translation and rotation operators \cite{king1973general,macphie1987rotational,dassios2012ellipsoidal} may be investigated in future work to accelerate this step to also be computed in $\mathcal{O}(Mp^3 \log p)$ work. 
\end{itemize}

Combining these results with our complexity analysis for far-field interactions, we can conclude that our methodology allows us to compute all-to-all interactions for a dense particulate system in $\mathcal{O}(N)$ work. The accelerations proposed above then serve to reduce the algorithmic constants for near and self-interactions. As demonstrated by our numerical experiments in Section \ref{sec:numerical}, given a target accuracy $\varepsilon$, our choice of $p$ only depends on $\varepsilon$ and particle eccentricities: for low to moderate aspect ratios, we observe rapid and uniform spectral decay of the relative evaluation errors; for aspect ratios around $8$ or higher, resolving the error on a boundary layer near the surface requires increasing $p$ considerably. 

Overall, this fast evaluation framework extends the benefits and complexity of the spectral, global basis function analysis approach of \cite{corona2018boundary}, which can be contrasted with panel-based approaches for axisymmetric surfaces \cite{young2012high} and slender bodies \cite{MALHOTRA2024112855}; for particles of low to moderate aspect ratio, we expect the method proposed here to require considerably less degrees of freedom to achieve a desired accuracy.  

\subsection{Matrix-vector evaluation algorithm}\label{ssc:matvec}

We now consider an example for the first task in our problem setup: evaluating the double layer potential at all particle surfaces. For each $\Gamma_i$, we calculate its contribution to the double layer potential evaluated at all target points in the system of $M$ spheroids. This routine is outlined in Algorithm \ref{alg:DL}: for each source spheroid, we compute normalized distances to all spheroids, separate them into far and near, and evaluate the resulting potential $u_i(\bm{x})$ by separating it into three pieces: $u_i^{far}$ (computed using smooth quadrature), $u_i^{self}$ (on-surface potential evaluation) and $u_i^{near}$ (off-surface potential evaluation). 

As mentioned previously, far-field computation (Lines $6-8$) can be accelerated by replacing direct quadrature evaluation with a Fast Multipole Method e.g. \cite{greengard2012fmmlib3d}. This would involve computing $u^{far}$ for all particles on a separate for loop, making sure to remove or omit self and near particle interactions from the FMM algorithm. 

\begin{algorithm}
\caption{Evaluate $u(\bm{x})=\DL[\mu](\bm{x})$}
\label{alg:DL}

\textbf{Inputs: } Densities $\{\mu_i\}_{i=1}^M$, Surfaces $\{\Gamma_i\}_{i=1}^M$, Order $p$, Parameter $\eta$ \\
\textbf{Output: } Layer potential $u(\bm{x})$

\begin{algorithmic}[1]
    \For{$i=1,\dots,M$}
    \State{ $\widehat{\mu_i} = \mathrm{Fast\_Spheroidal\_Transform(\mu_i)}$  (if not provided in inputs)}

    \State{Compute normalized distances $d_{ij} = d(\Gamma_i,\Gamma_j) / \mathrm{diam}(\Gamma_i)$}
    \State{Find far and near particle indices $J_i^{far} = \{j \ | \ d_{ij} \geq \eta\}$ and $J_i^{near} = \{j \ | \ d_{ij} < \eta\}$}
    \State{Evaluate $u_i(\bm{x})=\DL[\mu_i](\bm{x})$ in three regions:}
    \vskip 5 pt
    \For{$j \in J_i^{far}$}
        \State{Compute $u^{\text{far}}_i(\bm{x})$ using smooth quadrature at all $x \in \Gamma_j$.}
    \EndFor

    \vskip 5 pt
    \For{$j \in J_i^{near}$}
        \If{$j = i$}
            \State{Compute coefficients $\widehat{u_i} = \widehat{\DL\left[\mu_i\right]}$ \Comment{on-surface evaluation Sec \ref{sec:DL_eval}}}
            \State{$u^{\text{self}}_i(\bm{x}) = \mathrm{Inverse\_Fast\_Spheroidal\_Transform(}\widehat{u_i}\mathrm{)}$}
        \Else 
            \State{Compute solid harmonic coefficients of $\DL\left[\mu_i\right](\bm{x})$ 
            \State{Evaluate solid harmonic expansion $u^{\text{near}}_i(\bm{x})$} \Comment{off-surface evaluation Sec \ref{sec:DL_eval}}}    
        \EndIf
    \EndFor

    \vskip 5 pt
    \State{$u_i(\bm{x}) = u^{\text{far}}_i(\bm{x}) + u^{\text{near}}_i(\bm{x}) + u^{\text{self}}_i(\bm{x})$}
    \State{$u(\bm{x}) = u(\bm{x}) + u_i(\bm{x})$}
    \EndFor
    \State{\textbf{Return} $u(\bm{x})$}
\end{algorithmic}
\end{algorithm}

\subsubsection*{Spheroid distance calculation}

The normalized distance computation in Line $3$ of Algorithm \ref{alg:DL} is performed in two stages to ensure computational efficiency: First, we calculate lower bounds on distances using circumspheres for each spheroid. We then further refine the calculation for potentially nearby spheroids using the moving balls algorithm described in \cite{girault2022movingballs}. By definition, all targets points on a spheroid tagged as ``far'' are on the far-field, but of course, only a subset of points on a ``near'' particle are on the near-field. In our implementation, we have prioritized fast distance computation and processing of far and near computation at the particle level. Alternatively, exact distances of individual targets on ``near'' particles could be computed if processing all ``far'' points using the smooth rule is deemed advantageous. 

\subsubsection*{Other layer potentials and their derivatives}

A few notable differences arise when evaluating $\SL, \SL'$ and $\DL'$ instead of $\DL$ in Algorithm \ref{alg:DL}. Firstly, when evaluating $\SL$ in the ``near'' regime, the formulas to obtain $\SL[\sigma_i]$ depend on coefficients of $\sigma_i$ in the modified basis $\left\{\frac{Y_n^m}{\sqrt{u_0^2\pm v^2}}\right\}$. These coefficients $\widetilde{\sigma_n^m}$ can be obtained from a fast spheroidal transform applied to $\sigma_i \sqrt{u_0^2\pm v^2}$.  

As shown for $\SL'$ in Theorem \ref{thm:dSLthm}, on-surface formulas for $\SL'$ and $\DL'$ follow the same pattern as for layer potentials: we can obtain coefficients for a spheroidal or modified spheroidal basis from the coefficients for $\sigma_i$ via a diagonal operator, and then we can evaluate our operator on the surface using a fast inverse transform. However, for near off-surface evaluation (e.g. evaluation at targets on a neighboring spheroid), the normal vector at the target generally has components in $u,v$ and $\phi$ in spheroidal coordinates relative to the source spheroid. For this reason, evaluating  derivatives off-surface requires a formula to compute the gradient of solid spheroidal expansions, which is obtained for $\SL'$ in  Equation \eqref{eq:dSLevalformula}. 

\subsubsection*{Particle-to-target algorithm}

A modification of Algorithm \ref{alg:DL} can be made to evaluate a layer potential or its derivatives at a user-prescribed set of target points. Key differences are that the distance computation and classification between near and far is now done for each target point, and in general, on-surface evaluation (Line $12$) is done directly, unless there are enough target points on a given surface to warrant a different approach (e.g. fast transform and interpolation to targets). 

\subsection{Stokes potentials}
\label{sec:Stokes}

We now discuss the extension of the spheroidal-harmonic BIO evaluation scheme to Stokes layer potentials. Recall that, in the low Reynolds number limit, the velocity field $\bm{u}(\bm{x})$ of a Newtonian fluid satisfies the Stokes equations:
\begin{equation}
     -\mu\Delta \bm{u}+\nabla p=0, \quad \nabla\cdot\bm{u}=0 \quad \text{in} \quad \Omega,
     \label{eq:Stokes}
\end{equation}

where $p$ is the pressure, $\mu$ is the fluid viscosity and $\Omega$ denotes the fluid domain. Analogous to the Laplace case, Stokes layer potentials can be used to express the solution of \eqref{eq:Stokes} in terms of vector density functions defined on the domain boundary $\Gamma$. In a single-layer potential ansatz, for example, the fluid velocity and pressure are expressed as

\begin{align}
    \bm{u}(\bm{x}) = \SLt[\bm{\sigma}](\bm{x}) & := \int_\Gamma G_\text{stk} (\bm{x},\bm{y})\bm{\sigma}(\bm{y}) d\Gamma(\bm{y}) \quad \text{and} \\ p(\bm{x}) = \mathcal{P}[\bm{\sigma}](\bm{x}) & := \int_\Gamma P (\bm{x},\bm{y})\cdot\bm{\sigma}(\bm{y}) d\Gamma(\bm{y}),
\end{align}
%
%
%
where $\bm{\sigma}$ is an unknown vector density function defined on $\Gamma$ that needs to determined by imposing given boundary conditions, and the Stokes free-space Green's function and the pressure kernel are given by 
$$G_\text{stk}(\bm{x},\bm{y})=\frac{1}{8\pi\mu}\left(\frac{\bm{I}}{|\bm{r}|}+\frac{\bm{r}\otimes\bm{r}}{|\bm{r}|^3}\right) \quad \text{and}\quad  P (\bm{x},\bm{y}) = \frac{1}{4\pi}\frac{\bm{r}}{|\bm{r}|^3}, \quad \bm{r} = \bm{x}-\bm{y}.$$

%

When $\Gamma$ is a spheroid, the goal is to evaluate these Stokes BIOs using the formulas derived earlier. This is accomplished by simply invoking the standard result that Stokes BIOs can be written in terms of Laplace layer potentials and their derivatives \cite{tornberg2008fmm3dstokes}. For example, the single layer potential can be expressed as 
\begin{equation}
    \SLt[\bm{\sigma}](\bm{x})_k = \frac{1}{2}\left\{
    \SL [\sigma_k(\bm{y})](\bm{x})-
    \sum_{j=1}^3\left(
    x_j \ \bm{e}_k(\bm{x}) \cdot \nabla \SL [\sigma_j](\bm{x})\right)
    + \bm{e}_k(\bm{x}) \cdot \nabla \SL [\bm{y}\cdot\bm{\sigma}](\bm{x})
    \right\},\quad k=1,2,3,
    \label{eq:stokes}
\end{equation}
which implies that Stokes single layer evaluation can be reduced to three Laplace single layer evaluations and four gradient of single layer evaluations. For a Stokes single layer MatVec algorithm, we compute far-field interactions using smooth integration of the Stokeslet kernel, which can be accelerated via a Stokes FMM. We then incorporate near-field interactions using our formulas for Laplace and Eq \eqref{eq:stokes}. To evaluate interactions on the near-field for $\nabla \SL$, we need two components: 

\begin{itemize}
    \item \emph{Neighbor particle interactions.} We use the evaluation formula in Equation \eqref{eq:dSLevalformula}, replacing $\bm{\nu}(\bm{x})$ with $\bm{e}_k(\bm{x})$, expressed in spheroidal coordinates. This evaluates $\frac{\partial \SL}{\partial x_k}$ at neighbor particle targets.  
    \item \emph{Self-interactions.} We take the appropriate limit (exterior or interior) as $u \rightarrow u_0$ and apply it to Equation \eqref{eq:dSLevalformula}; this yields an evaluation formula for self-interactions for each partial derivative. 
\end{itemize}

Similarly, the pressure field can be evaluated from gradient of Laplace single layer potentials by writing it as
$$\mathcal{P}[\bm{\sigma}](\bm{x})=\sum_{j=1}^3\bm{e}_j(\bm{x})\cdot\nabla \mathcal{S}[\sigma_j](\bm{x}).$$

We note that similar formulas writing the Stokes double-layer $\DL_{stk}$ in terms of Laplace $\DL, \nabla \DL$ BIO evaluations as in \cite{bagge2021highly} may be readily implemented with this methodology.  

\subsubsection*{Example: Interfacial Stokes flow}
To test our spheroidal BIO evaluation scheme for Stokes flow, we consider an interfacial Stokes flow problem as illustrated in Figure \ref{fig:stokes flow}. In this setup, fluid-enclosing spheroidal particles (both oblate and prolate) are suspended in a viscous fluid confined by a prolate spheroid. The flow is driven by interfacial forces acting on the boundary of each of these particles, defined by $\bm{f}^\text{int} = [\![-p\bm{n} + (\nabla\bm{u} + \nabla\bm{u}^T )\cdot\bm{n}]\!]$, that is, the jump in hydrodynamic stresses in the normal direction. In this example, we set $\bm{f}^\text{int}(\bm{y}) = H(\bm{y})\bm{n}(\bm{y})$ for $\bm{y}\in\Gamma$, where $H$ is the mean curvature. 

We evaluate the curvature-induced flow inside the confinement at a given time using the Stokes Single Layer formulation: 
\[\bm{u}(\bm{x}) = \SLt[\bm{f}^\text{int}](\bm{x})\]
and plot the streamlines at a given $x_1x_2-$plane to observe how particles intersecting the plane impact the flow (Figure \ref{fig:stokes flow}(right)).


\section{Numerical results}
\label{sec:numerical}
In this section, we conduct a series of numerical experiments to test our proposed BIO evaluation framework. First, we test its performance in the solution of Laplace Dirichlet and Neumann boundary value problems via standard indirect integral equation formulations leading to Fredholm equations of the second kind. In both cases, the solution can be recovered with uniform, spectral accuracy at target points arbitrarily close to particle surfaces. We then conduct a stress test centered on the behavior of our framework when we increase particle aspect ratio and bring particles closer to each other. We present results from these tests focusing on the spatial distribution of the error and conditioning of commonly used integral equation operators.  

\subsection{Laplace boundary value problems on spheroidal suspensions}\label{ssc:LaplaceBVP}

A set of non-overlapping $M$ spheroids with boundaries $\Gamma_i$ is defined by shape information, i.e. prolate or oblate and parameters $u_0^i$ and $a_i$, centers $\bm{c}_i$ and rotation from a reference configuration encoded by a normalized quaternion $\bm{q}_i$. We denote the domain exterior to $\cup_{i=1}^M \Gamma_i$ as $\Omega$. We generate a harmonic function in $\Omega$ by placing $S$ point sources of strength $\tilde{\mu}_j$ located at points $\bm{x}_j$ inside these spheroids (point charges inside perfectly conducting spheroidal particles in an electrostatics context). The corresponding harmonic potential in $\Omega$ is then given by: 

$$f(\bm{x})=\sum_{j=1}^S \frac{\tilde{\mu}_j}{4\pi||\bm{x}-\bm{x}_j||},\quad x\in\Omega$$

We use instances of this general setup to test performance of our integral evaluation scheme in the solution of exterior Dirichlet and Neumann problems. In both cases, we

\begin{enumerate}
    \item Define an indirect integral formulation based on Laplace layer potentials for unknown integral densities $\sigma_i$. 
    \item Use the MatVec algorithm in Section \ref{ssc:matvec} for a range of values of spheroidal harmonic order $p$ to solve the corresponding discretized linear system to near machine precision, computing integral densities $\sigma_i$ at discretization points.
    \item Compare $f(x)$ with the integral operator evaluated at shells of target points set at distances approaching each $\Gamma_i$ exponentially. 
\end{enumerate} 

\begin{figure}
    \centering
    \includegraphics[width=\textwidth]{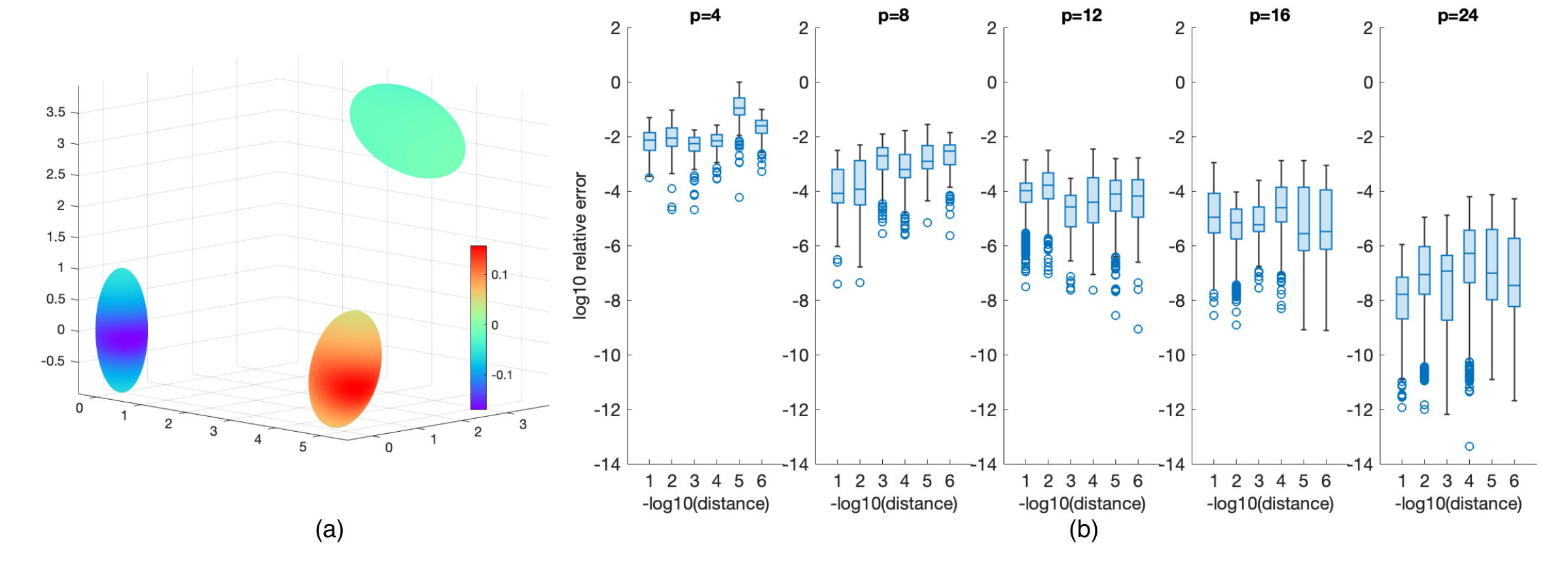}
    \caption{{\em \small (a) The system of prolate spheroids centered at $(0,0,0)$, $(5,0,0)$, $(3.2,3.2,3.2)$ with $u_0=1.1,1.2,1.3$ respectively, and with random orientations, on which the exterior Dirichlet problem is solved. Point sources of random strengths $-0.5\leq \tilde{\mu}_j\leq 0.5$ are placed at random near the center of each spheroid, $j=1,...,6$, generating a potential on the surface of each spheroid that is plotted by the heatmap in this figure. This potential is used as Dirichlet boundary data.
    (b) We solve the BIE \eqref{testprob} with equation \eqref{completion} as a completion term, at distance $10^{-k}$ away from the surface of each spheroid, $k=1,...,6$. The solution for this BIE is recovered with relative error that decays spectrally with order $p$; this is maintained as evaluation points approach the surface and the integral kernel becomes singular.
    }}
    \label{fig:spheroids}
\end{figure}

\subsubsection*{Laplace exterior Dirichlet problem}

We consider a system of $3$ prolate spheroids with aspect ratio around $1.8$
with different orientations,
and randomly generate $2$ point sources inside each spheroid, with strengths between $[-0.5,0.5]$, for a total $S=6$. The resulting potential on particle surfaces is plotted in Figure \ref{fig:spheroids}(a). We then solve the exterior Dirichlet problem
 $$\begin{cases}
     \Delta u(\bm{x})=0 & \bm{x}\in\Omega\\
     u(\bm{x})=f(\bm{x}) & \bm{x}\in\bigcup_i\Gamma_i
 \end{cases}$$
using the \emph{completed Double Layer} formulation  $u(\bm{x})= \DL^+[\sigma](\bm{x}) +\mathcal{C}[\sigma](\bm{x})$, where $C$ denotes an additional term that completes the null space of $\frac{1}{2}I+\DL$, making the resulting BIE uniquely solvable:
\begin{equation}\label{testprob}
    \frac{1}{2}\sigma(\bm{x}) + \DL[\sigma](\bm{x}) + \mathcal{C}[\sigma](\bm{x}) = f(\bm{x}), \quad \bm{x}\in \bigcup_i \Gamma_i
\end{equation}
There are a number of options for $C$; a widely used and arguably the simplest one corresponds to placing point sources of strengths equal to the integral of $\sigma_i$ at each spheroidal center, that is, 

\begin{equation}\label{completion}
\mathcal{C}_I[\sigma](\bm{x}) = \sum_{i=1}^M\frac{1}{||\bm{x}-\bm{c}_i||}\int_{\Gamma_i} \sigma_i(\bm{y}) \rmd S
\end{equation}

We use our evaluation framework to test how the choice of completion terms affects conditioning and GMRES iteration counts as a function of aspect ratio and inter-particle distance in Section \ref{ssc:StressTest}.

\subsubsection*{Laplace exterior Neumann problem}

We then consider the exterior Neumann boundary value problem
given the same potential $f(x)$ in $\Omega$. For this test case, we evaluated performance on a mixed system of two oblate and one prolate spheroids, whose configuration and surface potential are plotted in Figure \ref{fig:fig3}(a). From an electrostatics interpretation, the point charges create an electric field that passes through the surface of each spheroid, 

$$\bm{E}(\bm{x}) = -\nabla f(\bm{x}) = -\sum_{j=1}^6\frac{\tilde{\mu}_j}{4\pi||\bm{x}-\bm{x}_j||^2}\frac{\bm{x}-\bm{x}_j}{||\bm{x}-\bm{x}_j||}$$

And so, the boundary data corresponds to setting outgoing fluxes in $\Gamma_i$, that is, 

$$\begin{cases}
    \Delta u(\bm{x})=0 & \bm{x}\in\Omega\\
    \frac{\partial u}{\partial n_+}(\bm{x}) =-\bm{E}(\bm{x})\cdot \bm{n}(\bm{x}) & \bm{x}\in\bigcup_i\Gamma_i
\end{cases}$$

We solve this using the single layer formulation 

\begin{equation}
u(\bm{x})=\SL^+[\rho](\bm{x})\label{eq:neusol}
\end{equation}

leading to the boundary integral equation 


$$-\frac{1}{2}\rho(\bm{x})+\mathcal{S^{\prime}}[\rho](\bm{x})=-\bm{E}(\bm{x})\cdot \bm{n}(\bm{x}),\quad \bm{x}\in\bigcup_i\Gamma_i$$

\subsubsection*{Numerical results}

\begin{figure}
    \centering
    \includegraphics[width=\textwidth]{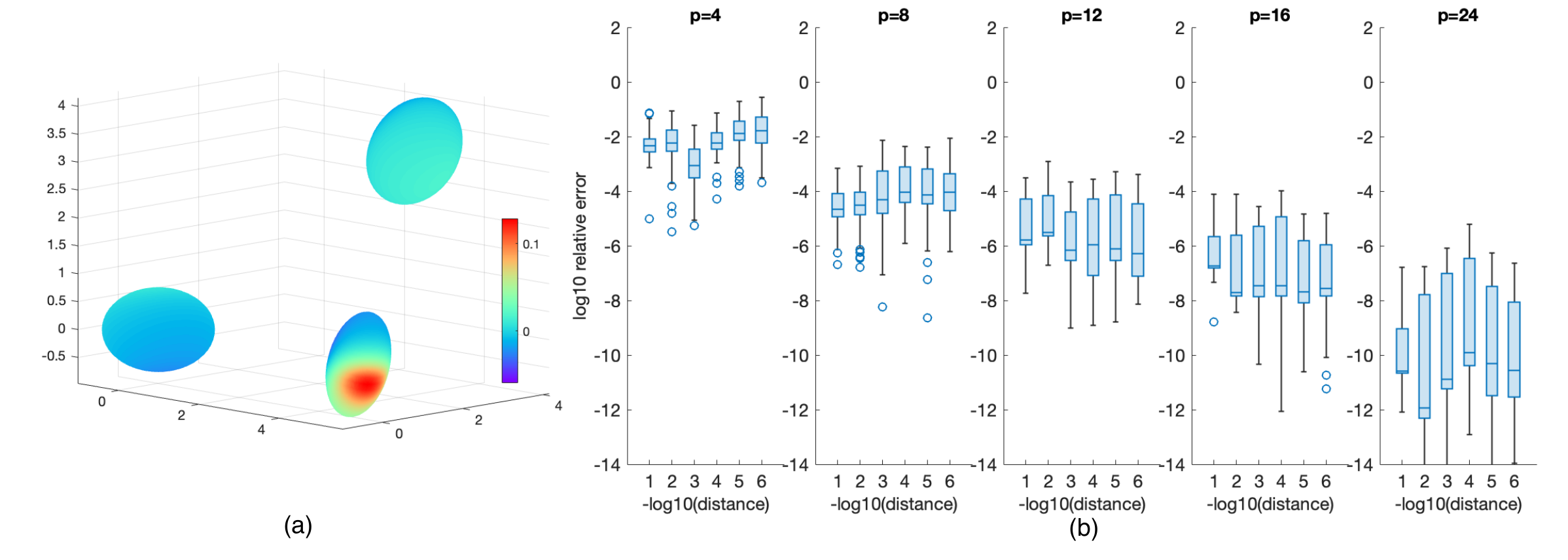}
    \caption{{\em \small We test the accuracy of the evaluation method on a system of two oblate spheroids centered at $(0,0,0)$ ($u_0=1.1$), $(3.2,3.2,3.2)$ ($u_0=1.3$), and one prolate spheroid, centered at $(5,0,0)$ ($u_0=1.2$). Two point charges are randomly generated inside each spheroid as in Figure \ref{fig:spheroids}; the resulting potential on the surface of each spheroid is plotted in (a). We solve the exterior Neumann boundary value problem using ansatz \ref{eq:neusol}, at distance $10^{-k}$ away from the surface of each spheroid, $k=1,...,6$. At each distance, the method demonstrates spectral decay of relative error as the order $p$ increases, shown in (b).   
    }}
    \label{fig:fig3}
\end{figure}

The results in Figure \ref{fig:spheroids}(b) and Figure \ref{fig:fig3}(b) are representative of the performance of our evaluation framework across a range of examples involving integral operator evaluation in suspensions of prolate and oblate spheroids: for suspensions with moderate aspect ratios, the solution is recovered with fairly uniform, spectral rate regardless of target point proximity to a spheroidal surface. It is important to note that this showcases performance of our method both in solving the BIE to high accuracy and in evaluating the solution on a set of target points, regardless of how close they are to spheroidal surfaces. 

\subsection{Stress tests on spheroidal evaluation framework}\label{ssc:StressTest}


In this section, we conduct a series of stress tests to determine how the accuracy of the method as well as the condition number of the integral operator behave when one or more particles' aspect ratio increases (prolates become more elongated, while oblates become increasingly squat), especially in cases where two or more particles come close to each other. The aspect ratio of a spheroid is the major axis over minor axis. For a prolate spheroid, it is $R=C/A=u_0/\sqrt{u_0^2-1}$, and as $R$ increases the spheroid becomes elongated along its major axis and approaches a line segment (similar to a slender rod). For an oblate spheroid, the aspect ratio is $R=\sqrt{u_0^2+1}/u_0$, and the particle ``flattens", approaching a disk.

\subsubsection*{Accuracy}

\begin{figure}
        \centering
        \hspace{-.2in}\includegraphics[width=1.04\textwidth]{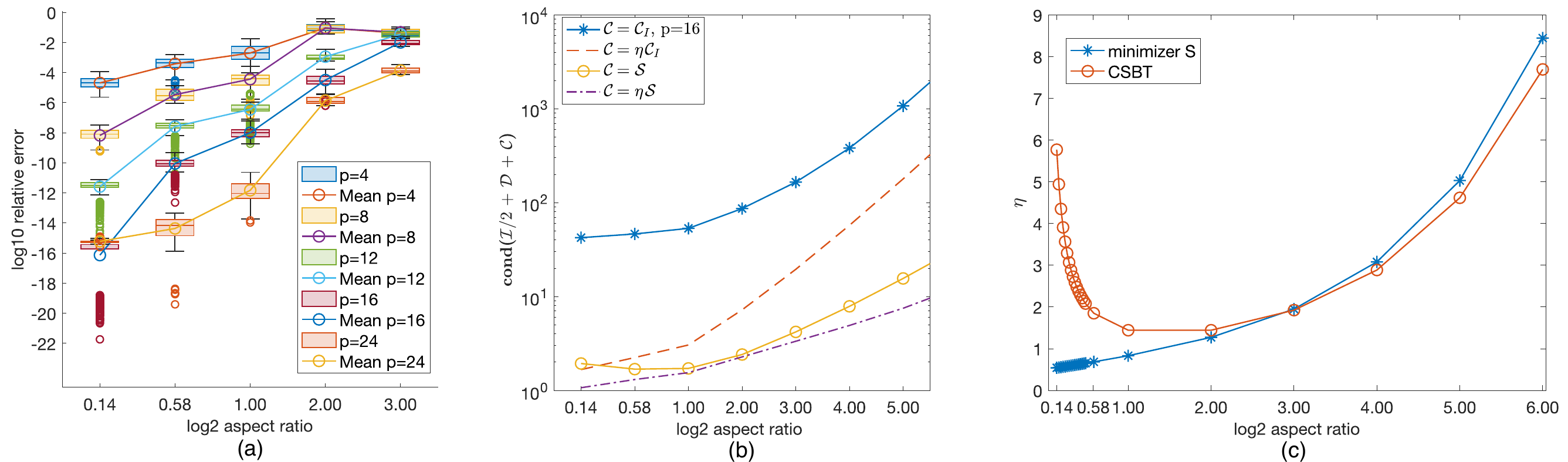}
        \caption{{\em \small (a) We solve an exterior Dirichlet problem with one prolate particle containing two point charges near its center, and compare the BIE solution to the true potential of the system on a order 16 Gauss-Legendre grid $0.5$ away from the surface of the spheroid. The completion term \eqref{completion} is used. Relative error of the spectral method decays exponentially with $p$ until aspect ratio reaches $R=4$, at which point the error plateaus at three digits. (b) Condition number of the exterior Dirichlet self-evaluation operator on one prolate spheroid, using the completion terms $\mathcal{C}_I$, $\eta\mathcal{C}_I[\sigma]$, $\SL[\sigma]$, and $\eta\SL[\sigma]$. $\eta$ is taken to be the minimizer of the condition of this operator with positivity constraint. The rank 1 completion $\mathcal{C}_I$ and $\eta\mathcal{C}_I$ quickly become ill-conditioned as aspect ratio increases despite condition-minimization, while both $\SL$ and $\eta\SL$ remain well-conditioned. (c) The minimizing $\eta$ for $\SL$ is compared to scalar $1/(2\varepsilon\log(\varepsilon^{-1}))$ from CSBT, where $\varepsilon=1/R$. For $R\geq 4$, the two closely agree.}} \vspace{-0.1in}
        \label{fig:fig4}
\end{figure}

 To observe the extent to which our proposed method handles high-aspect-ratio particles, we first limit the system to one prolate spheroid and solve the exterior Dirichlet problem in Section \ref{ssc:LaplaceBVP}. We again place two point sources inside with randomized location and strength, and observe the error at a shell of target points placed at a distance of $0.5$ away from the particle surface for a range of increasing aspect ratios, going from $\simeq 1$ to $8$. 
 We observe in Figure \ref{fig:fig4}(a) that our evaluation method displays spectral convergence up until an aspect ratio $R=4$. However, as we further increase the aspect ratio, the error at these specific target points plateaus. The same happens for an oblate spheroid. 

To further investigate this behavior and check the spacial distribution of error, we place four identical upright prolate spheroids in a two-by-two lattice, and plot the log error on the $x$-$z$ plane when we increase their aspect ratio. 
To add complexity of the problem, we placed 21 randomly generated charges inside each particle. We first notice that for each fixed aspect ratio, error decreases spectrally overall as the order of our method increases. However, for aspect ratios $R \geq 4$, a neighborhood close to the surface seems to retain high error that requires much larger p to resolve. Separate tests ran on oblate and mixed spheroid suspensions show similar error behavior. 

\begin{figure}
    \centering
    \includegraphics[width=\linewidth]{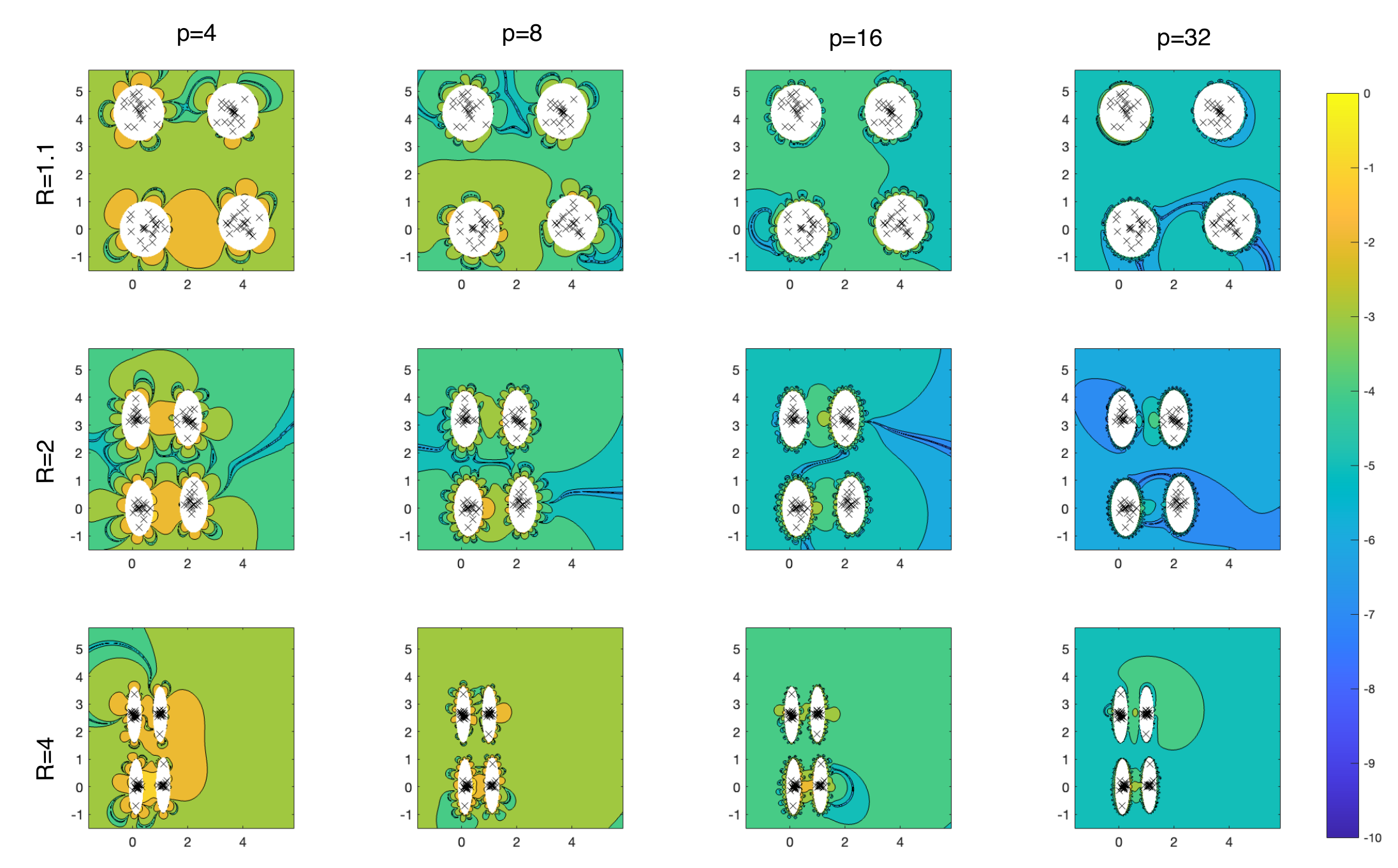}
    \caption{{\em \small Error of potential evaluation is plotted for a two-by-two lattice of four identical prolate spheroids, each with 21 random charges inside. More elongated spheroids are brought closer together to maintain proximity of evaluation. At fixed aspect ratio, the method is spectrally convergent overall. As aspect ratio increases, larger error persist close to the particles, preventing the overall error from decreasing, hence the lack of convergence seen in Figure \ref{fig:fig4}(a).}}
    \label{fig:err_contour}
\end{figure}

\subsubsection*{Condition number}

In the context of completed double layer formulations for Laplace and Stokes Dirichlet problems, it is well known that the conditioning of the operator using $\mathcal{C}_I$ can considerably deteriorate for elongated prolates and slender bodies, as can be seen in Figure \ref{fig:fig4}(b); this results in higher iteration counts for Krylov subspace methods like GMRES. A commonly used alternative is $\mathcal{C}=\SL$; a scaled $\SL$ with the scaling factor based on relative cross section radius is investigated in \cite{MALHOTRA2024112855} for slender bodies, where the method is shown to converge as aspect ratio increases. Other alternatives include a generalization of $\mathcal{C}_I$ (presented for the Stokes case in \cite{keaveny2011applying}) using a segment of sources along the object's centerline. 

In this section, we use our evaluation framework to compare unscaled and scaled versions of $\mathcal{C}_I$ and $\SL$; for scaled completion terms, we use an optimization solver to numerically find the $\eta$ that minimizes condition number of the completed Laplace operator for one particle. In Figure \ref{fig:fig4}(b) we show that despite the minimizing process, the operator with $\mathcal{C}_I$ still becomes much more ill-conditioned than the double layer completed with either unscaled or scaled $\SL$, as the aspect ratio of the spheroid increases. 

This work motivates a heuristic scaling scheme, which we summarize in Table \ref{tab:eta_formulas}. If $\mathcal{C}_I$ is used, we find that for prolates and oblates of low or moderate aspect ratio, scaling it by a factor proportional to the inverse of the surface area, $\varepsilon/A(\varepsilon)$, achieves minimal condition number for $\mathcal{K}$, reducing conditioning by about one order of magnitude. For high aspect ratio oblates, the minimizer approaches a constant of about $0.06$. If, on the other hand, $\SL$ is used, we find that the slender body limit (SBLT) scaling $\eta(\varepsilon) = \frac{1}{2 \varepsilon \log(1/\varepsilon)}$ as in \cite{MALHOTRA2024112855} closely matches the minimizer for elongated prolates with $R \geq 3$, while $0.5$ is a good match for prolates with $R \leq 3$. The scaling factor for $\SL$ on the oblate case approaches $\eta=1$ for $R>3$ (as the oblate approaches a flat disk); for simplicity, we assign $\eta=1$ for $R>3$ and vary the scalar linearly for $R\leq 3$, since the condition number is small and scaling does not impact it drastically.

\begin{center}
   \begin{tabular}{ c | c | c | c }
   {\em Completion term}&{\em Aspect ratio} & {\em \textbf{Prolate}} & {\em \textbf{Oblate}} \\
   \hline 
   & & \tabularnewline[-0.5em]
  \multirow{ 2}{*}{$\mathcal{C}_I$}  & $R=\frac{1}{\varepsilon}\leq 3$ &
\multirow{ 2}{*}{$\frac{\varepsilon}{A^{\mathbb{P}}(\varepsilon)}$} 
     & $\frac{\varepsilon}{A^{\mathbb{O}}(\varepsilon)}$
     \\
     & &  \tabularnewline[-0.5em]
     & $R>3$ & & $0.06$\\
    & &  \tabularnewline[-0.5em]
   \hline
   & & \tabularnewline[-0.5em]
  \multirow{ 2}{*}{$\SL$} & $R\leq 3$ & $0.5$ & $\frac{1}{4\varepsilon}+\frac{1}{4}$\\
   & & \tabularnewline[-0.5em]
   & $R>3$ & $\frac{1}{2 \varepsilon \log(1/\varepsilon)}$ & $1$
   \end{tabular}
   \captionof{table}{{\em \small Value of $\eta$ as a function of the elongation $\varepsilon=\frac{1}{R}$ heuristically chosen to match the minimizer of condition number of integral operator $\mathcal{K}=\frac{1}{2}\mathcal{I}+\DL + \mathcal{C}$. 
   Scaling for completion flow $\mathcal{C}_I$ relates to surface area of the particle for both prolate and oblate spheroids; scaling for $\SL$ follows CSBT for elongated prolate spheroids, and approaches 1 for high aspect ratio oblates.}}
   \label{tab:eta_formulas}
\end{center} \vspace{0.05in} 

We demonstrate the impact of unscaled and scaled completion terms as described above in terms of how condition number affects GMRES iteration counts. For this purpose, we consider the lattice system of four identical prolate spheroids in the previous section. For each completion $\mathcal{C}=\mathcal{C}_I,\ \eta\mathcal{C}_I,\ \SL,\ \eta\SL$, the number of GMRES iterations to solve $K\bm{\sigma}=f$ (the discretized matrix equation for the completed double layer integral equation $\mathcal{K}[\sigma]=f$) is reported in Table \ref{GMRES} for increasing aspect ratios. We proportionally decrease the distance between spheroids ($d$) based on the aspect ratio. For GMRES, the tolerance is set to 1e-10 and max iteration to 1e6. Note that the GMRES iteration count is upper bounded by the size of the matrix $K$, so the MATLAB algorithm takes the size of the matrix as max iteration when this parameter is set to a larger value. 

\begin{center}
    \begin{tabular}{c|c|c|c|c|c|c|c|c|c|c|c|c}
    &\multicolumn{3}{c|}{$\mathcal{C}_I$} & \multicolumn{3}{c|}{$\eta\mathcal{C}_I$} & \multicolumn{3}{c|}{$\SL$} & \multicolumn{3}{c}{$\eta \SL$}  \\
    \hline 
    $R\times d$ & $2$ & $1$ & $0.01$ & $2$ &  $1$ & $0.01$ & $2$ & $1$ & $0.01$ & $2$ & $1$ & $0.01$ \\
    \hline 
    $R=1.1$ &14 & 15 & 36 &13 & 14 & 33 & 10 & 11& 33 &10 & 12 & 34\\
    $R = 2$ &19 & 19 & 39 &17 & 17 & 35 & 12& 13 & 35 &14 & 16 & 37\\
    $R = 4$ &27 & 28 & 49 &26 & 27 & 51&16 & 18 & 39 & 15& 17 & 38\\
    $R = 8$ &44 & 45 & 65 &42 & 45 & 68 &23 & 25 & 51 &20 & 21 & 45\\
    $R = 16$ & 69 &75  &107 & 67 &72 &115 & 31 & 34 &68 & 24 &26 &52 \\
    $R = 32$ &103 & 109 & 170 &100 & 107 & 169 &40 & 45 & 91 &31 & 32 &  60\\
    $R = 64$ &131 & 150 & 223 &133 & 154 & 228 &47 & 57 & 124 &38 & 40 & 67
    \end{tabular} 
    \captionof{table}{{\em \small Number of iterations for MATLAB GMRES to converge with 1e-10 tolerance, to solve $K\bm{\sigma} = f$, for configurations with increasing aspect ratio and decreasing separation. $\mathcal{K}=\mathcal{I}/2+\DL+\mathcal{C}$, with $\mathcal{C}=\mathcal{C}_I,\ \eta\mathcal{C}_I,\ \SL,\ \eta\SL$, for methods of order $p=16$ and prolate spheroid aspect ratios of $1.1,\ 2,\ 4,\ 8,\ 16,\ 32,$ and $64$. The benefit in computation cost for $\eta \SL$ over $\mathcal{C}_I$ is evident especially for high elongation and close proximity.}}
    \label{GMRES}
\end{center}


For low to moderate aspect ratios, the results suggest the differences between the four completion terms are slight; iteration counts are a bit smaller for the unscaled or scaled single layer. For $R >4$, especially for cases in which the particles come closer together, we start to see the advantage of using $\eta\SL$. The number of iterations needed to converge for a solution $\sigma$ increases dramatically for completion operators $\mathcal{C}_I$ and $\eta\mathcal{C}_I$, somewhat significantly for $\SL$, but barely noticeably for $\eta\SL$. Adding scaling to completion $\mathcal{C}_I$, although decreasing the condition number by a magnitude, does not seem to decrease the GMRES count significantly; on the other hand, the scaled $\eta\SL$ gives much more stable and consistently lower iteration counts than $\SL$. 
For densely packed systems of elongated prolates, using a scaled $\SL$ to complete the double layer formulation is therefore highly favorable due to its controlled condition number growth. 

Running the same tests for oblate particles reveals that either unscaled and scaled $\SL$ (which become the same at high aspect ratio) yield favorable conditioning and result in the lowest iteration counts for the completed double layer. Our observations for $\mathcal{C}_I$ match those for prolates: scaling this term improves condition number significantly, but gains in terms of iteration counts are small.  





\section{Conclusions}
\label{sec:conclusions}
We have presented a fast, spectrally accurate method for the evaluation of BIOs on suspensions of prolate and oblate spheroids. The formulas derived in Section 3 allow for fast computation arbitrarily close to the surface of spheroids, resolving the high computational cost typical of near-singular quadrature methods. This, combined with FMM-accelerated evaluation of far-field interactions yields fast and accurate computation of scalar fields in tightly packed spheroidal suspensions. Through a number of experiments, we demonstrate reliable performance of our method in the solution of Laplace boundary value problems as well as evaluation of Laplace and Stokes potentials at arbitrary target locations. Our stress tests show practical spectral convergence up to moderately high aspect ratios. For Laplace Dirichlet problems, we use our framework to propose a completion term that greatly reduces conditioning and iteration counts across a wide range of aspect ratios for both prolates and oblates. 

The work contributed here is part of an ongoing research agenda addressing key computational challenges in boundary integral methods for particulate suspensions simulation with the aim to enhance our  ability to model problems of interest in soft matter systems. 
By extending the methodology in \cite{corona2018boundary,yan2019scalable}, we are now able to deploy fast, spectral techniques to evaluate layer potentials resulting from dense, polydisperse suspensions of spherical, prolate and oblate particles, while utilizing complementary techniques such as \cite{veerapaneni2011fast,broms2024method,nitsche2025corrected,bagge2021highly} to address non-spheroidal or non rigid particles and geometries. 

Ongoing and future work furthering this effort features the efficient implementation of state-of-the-art optimization solvers to a complementarity based approach for spheroidal collision resolution following the work in \cite{yan2019scalable} and the development of novel vector spheroidal harmonic bases and corresponding formulas for Stokes layer potential evaluation, as well as formulas for other kernels of interest (e.g. Screened Laplace as in \cite{kohl2023fast}). Ultimately, this work should enable a wide range of robust and reliable direct numerical simulation of systems of interest across application fields.

\section{Acknowledgements} 
We acknowledge support from NSF under grants DMS-2012424 and DMS-2309661.

\bibliographystyle{plainnat}
\bibliography{references}
\end{document}